\documentclass[a4paper, leqno, 12pt]{article}
\usepackage{amsmath,amsthm}
\usepackage{amssymb,latexsym}
\usepackage{enumerate}
\usepackage{color}
\usepackage{array}
\usepackage[all]{xypic}

\usepackage[utf8]{inputenc}

\usepackage[left=3cm,right=3cm,top=3.0cm,bottom=3.0cm]{geometry}

\newtheorem{thm}{Theorem}[section]
\newtheorem{cor}[thm]{Corollary}
\newtheorem{lm}[thm]{Lemma}
\newtheorem{prop}[thm]{Proposition}

\newtheorem{con}[thm]{Conjecture}

\theoremstyle{definition}

\newtheorem{rem}[thm]{Remark}

\numberwithin{equation}{section}

\font\Mssym=msbm9  scaled \magstep 1

\def\ZZ{{{\mbox{\Mssym Z}}}}

\def\NN{{\mathbb{N}}}

\def\CA{{\cal A}}
\def\CC{{\cal C}}

\def\CD{{\cal D}}

\def\CF{{\cal F}}
\def\CG{{\cal G}}

\def\CS{{\cal S}}

\def\CX{{\cal X}}

\def\epv {{$\mbox{}$\hfill ${\Box}$\vspace*{1.5ex} }}

\def\Im{\textnormal{Im}}

\def\snull{\textnormal{ }}

\def\ra{\rightarrow}

\def\ov{\overline}

\def\mod{\mbox{{\rm mod}}}

\def\KG{\textnormal{KG}}

\def\ind{\mbox{{\rm ind}}}

\def\Mod{\mbox{{\rm Mod}}}
\def\Ker{\mbox{{\rm Ker}}}
\def\End{\mbox{{\rm End}}}
\def\Hom{\mbox{{\rm Hom}}}

\def\Coker{\mbox{{\rm Coker}}}  \def\Coim{\mbox{{\rm Coim}}}
\def\supp{\mbox{\rm supp}}

 \def\Ind{{\rm Ind}} \def\ind{{\rm ind}} \def\ob{{\rm ob}}

\def\wh#1{\widehat{#1}}

\def\rad{\textnormal{rad}} 
 \def\Add{\textnormal{Add}} 

\def\ov#1{\overline{#1}} \def\snull{\textnormal{ }}
\def\und#1{\underline{#1}}

\def\op{\textnormal{op}}

\begin{document}

\baselineskip=17pt


\title{On Krull-Gabriel dimension and Galois coverings}

\author{Grzegorz Pastuszak${}^{*}$}

\date{}

\maketitle

\begin{center}Dedicated to the memory of Gena Puninski\end{center}

\renewcommand{\thefootnote}{}

\footnote{${}^{*}$Faculty of Mathematics and Computer Science, Nicolaus Copernicus University, Chopina 12/18, 87-100 Toru\'n, Poland, past@mat.uni.torun.pl.}

\renewcommand{\thefootnote}{\arabic{footnote}}
\setcounter{footnote}{0}


\begin{abstract}
Assume that $K$ is an algebraically closed field, $R$ a locally support-finite locally bounded $K$-category, $G$ a torsion-free admissible group of $K$-linear automorphisms of $R$ and $A=R\slash G$. We show that the Krull-Gabriel dimension $\KG(R)$ of $R$ equals the Krull-Gabriel dimension $\KG(A)$ of $A$. Furthermore, we determine the Krull-Gabriel dimension of locally support-finite repetitive $K$-categories. These results are applied to determine the Krull-Gabriel dimension of standard selfinjective algebras of polynomial growth. Finally, we show that there are no super-decomposable pure-injective modules over standard selfinjective algebras of domestic type.
\end{abstract}

\section{Introduction and main results}

Assume that $K$ is an algebraically closed field and $A$ a finite dimensional associative basic $K$-algebra with a unit. We denote by $\mod(K)$ the category of all finite dimensional $K$-vector spaces and by $\mod(A)$ the category of all finitely generated right $A$-modules. Let $\CF(A)$ be the category of all finitely presented contravariant $K$-linear functors from $\mod(A)$ to $\mod(K)$. It is a hard problem to describe the category $\CF(A)$ even if the category $\mod(A)$ is well understood. A natural approach to study the structure of $\CF(A)$ is via the associated \textit{Krull-Gabriel filtration} \cite{Po} $$\CF(A)_{-1}\subseteq\CF(A)_{0}\subseteq\CF(A)_{1}\subseteq\hdots\subseteq\CF(A)_{\alpha}\subseteq\CF(A)_{\alpha+1}\subseteq\hdots$$ of $\CF(A)$ by Serre subcategories, defined recursively as follows: 
\begin{enumerate}[\rm(1)]
	\item $\CF(A)_{-1}=0$,
	\item $\CF(A)_{\alpha+1}$ is the Serre subcategory of $\CF(A)$ formed by all functors having finite length in the quotient category $\CF(A)\slash\CF(A)_{\alpha}$, for any ordinal number $\alpha$,
	\item $\CF(A)_{\beta}=\bigcup_{\alpha<\beta}\CF(A)_{\alpha}$, for any limit ordinal $\beta$.
\end{enumerate} 
Following \cite{Ge1}, \cite{Ge2}, the \textit{Krull-Gabriel dimension} $\KG(A)$ of $A$ is the smallest ordinal number $\alpha$ such that $\CF(A)_{\alpha}=\CF(A)$, if such a number exists, and $\KG(A)=\infty$ if it is not the case. If $\KG(A)=n$ for some natural number $n$, then the Krull-Gabriel dimension of $A$ is \textit{finite}. If $\KG(A)=\infty$, then the Krull-Gabriel dimension of $A$ is \textit{undefined}. 

The interest in the Krull-Gabriel dimension $\KG(A)$ of an algebra $A$ has at least three motivations. The general one is that the Krull-Gabriel filtration of the category $\CF(A)$ leads to a hierarchy of exact sequences in $\mod(A)$ where the almost split sequences form the lowest level, see \cite{Ge1}. The second one is that there exists a strong relation between $\KG(A)$ and the transfinite powers $\rad_{A}^{\alpha}$ of the radical $\rad_{A}$ of $\mod(A)$, see \cite{Kr} and \cite{Sch1,Sch2,Sch3}. For example, H. Krause proves in \cite[Corollary 8.14]{Kr} that if $\KG(A)=n\in\NN$, then $\rad_{A}^{\omega(n+1)}=0$ where $\omega$ is the first infinite cardinal. Moreover, J. Schr\"oer conjectures in \cite{Sch3} that $\KG(A)=n\geq 2$ if and only if $\rad_{A}^{\omega(n-1)}\neq 0$ and $\rad_{A}^{\omega n}=0$. This conjecture is confirmed for several important classes of algebras, see for example Section 1 in \cite{Sk5} for the list. 

The third motivation to study the Krull-Gabriel dimension of an algebra is the following conjecture due to M. Prest, see \cite{Pr}, \cite{Pr2}.

\begin{con} A finite dimensional algebra $A$ is of domestic representation type if and only if the Krull-Gabriel dimension $\KG(A)$ of $A$ is finite.  
\end{con}

We refer to Chapter XIX of \cite{SiSk3} for the precise definitions of finite, tame and wild representation type of an algebra, as well as the stratification of tame representation type into domestic, polynomial and nonpolynomial growth (introduced in \cite{SkBC}). 

We recall some important results on the Krull-Gabriel dimension of an algebra. All of them support the above conjecture. Indeed, M. Auslander proves in \cite[Corollary 3.14]{Au} that $\KG(A)=0$ if and only if the algebra $A$ is of finite representation type. H. Krause shows in \cite[11.4]{Kr2} that $\KG(A)\neq 1$ for any algebra $A$. W. Geigle proves in \cite[4.3]{Ge2} that if $A$ is a tame hereditary algebra, then $\KG(A)=2$. A. Skowro\'nski shows in \cite[Theorem 1.2]{Sk5} that if $A$ is a cycle-finite algebra \cite{AsSk1}, \cite{AsSk2} of domestic representation type, then $\KG(A)=2$, see also \cite{Mal}. M. Wenderlich proves in \cite{We} that if $A$ is a strongly simply connected algebra \cite{Sk1}, then $A$ is of domestic type if and only if $\KG(A)$ is finite. R. Laking, M. Prest and G. Puninski prove in \cite{LaPrPu} a renowned result that string algebras \cite{SkWa}, \cite{BuRi} of domestic representation type have finite Krull-Gabriel dimension, see also \cite{Pu3}, \cite{PrPu}. 

There exist algebras with Krull-Gabriel dimension undefined. Indeed, the results of \cite[Chapter 13]{Pr} and \cite{GP} imply that the Krull-Gabriel dimension of strictly wild algebras and wild algebras, respectively, is undefined. J. Schr\"oer proves in \cite[Proposition 2]{Sch2} that this is also the case for nondomestic string algebras. W. Geigle shows in \cite{Ge1}, see also \cite{Ge2}, that tubular algebras \cite{Ri} have no Krull-Gabriel dimension. This also holds for pg-critical algebras, strongly simply connected algebras of nonpolynomial growth (see \cite{KaPa} and \cite{KaPa2}) and some algebras with strongly simply connected Galois coverings (see \cite{KaPa3}). We note that the results of \cite{KaPa,KaPa2,KaPa3} show, among other things, that the width of the lattice of all pointed modules \cite{PPT}, \cite{KaPa2} over any of these algebras is infinite. This yields that their Krull-Gabriel dimension is undefined, see for example \cite{Pr2}.

The Krull-Gabriel dimension is defined for any locally bounded $K$-category. The definition is similar to the one presented above, see Section 2 and Section 4 for the details. In fact, the Krull-Gabriel dimension is defined for any skeletally small abelian category, see Section 3. For example, in \cite{BobKr} the authors determine the Krull-Gabriel dimension of discrete derived categories. 

This paper is devoted to give a proof and derive some consequences of the following theorem.

\begin{thm}
Assume that $R$ is a locally support-finite locally bounded $K$-category and $G$ is a torsion-free admissible group of $K$-linear automorphisms of $R$. Assume that $A=R\slash G$ is the orbit category and $F:R\ra A$ the associated Galois covering. In this case we have $\KG(R)=\KG(A)$.
\end{thm}

Theorem 1.2 is a part of the more detailed Theorem 6.3. In Theorem 7.3 we determine the Krull-Gabriel dimension of locally support-finite repetitive $K$-categories. Specifically, it belongs to the set $\{0,2,\infty\}$. These two results are applied in Section 8 to determine the Krull-Gabriel dimension of standard selfinjective algebras of polynomial growth. Indeed, we prove in Theorem 8.1 that if $A$ is such an algebra, then $\KG(A)=2$ if $A$ is domestic and $\KG(A)=\infty$ otherwise. Note that this theorem supports Conjecture 1.1. As an application of Theorem 8.1 (and some results of \cite{KerSk}) we show in Corollary 8.2 that if $A$ is a standard representation-infinite selfinjective algebra of polynomial growth, then $\KG(A)=2$ if and only if the infinite radical $\rad_{A}^{\omega}$ is nilpotent. 

Although Theorem 6.3 is the main result of the paper, we view Theorem 7.3 and Theorem 8.1 as equally significant. 


Theorem 1.2 can be applied in some other situations as well. For example, A. Skowro\'nski studies in \cite{Sk3} algebras having strongly simply connected Galois coverings. Assume that $A$ is such an algebra, that is, $A=R\slash G$ where $R$ is a strongly simply connected locally bounded $K$-category and $G$ an admissible group of $K$-linear automorphisms of $R$. If $R$ is locally support-finite and $G$ is torsion-free, then Theorem 1.2 implies that $\KG(A)$ is finite if and only if $\KG(R)$ is finite. We conjecture this yields $\KG(A)$ is finite if and only if $A$ is of domestic type. This problem is left for further research.

The paper is organized as follows. In Section 2 we recall some information on Galois coverings of locally bounded $K$-categories and the associated pull-up and push-down functors. We conclude the section with Theorem 2.1 stating that if $R$ is a locally support-finite locally bounded $K$-category, $G$ an admissible torsion-free group of $K$-linear automorphisms of $R$ and $F:R\ra A\cong R\slash G$ the associated Galois covering, then the push-down functor $F_{\lambda}:\mod(R)\ra\mod(A)$ is a Galois covering of module categories, see Section 2 for definitions. This theorem is fundamental in proofs of our main results. Section 2 is based on \cite{Ga}, \cite{BoGa} and \cite{DoSk}.

Section 3 is devoted to the presentation of some facts about Serre subcategories of abelian categories, quotient categories and finite length objects (see especially Lemma 3.1). We define the Krull-Gabriel dimension of a skeletally small abelian category $\CC$. Among other things, we give in Lemma 3.2 a sufficient condition for an exact functor to preserve the Krull-Gabriel dimension. We apply the lemma in Section 6. Section 3 is based on \cite{Po}, \cite{McL} and \cite{Ge2}. 

In Section 4 we define the category $\CF(R)$ of all finitely presented contravariant $K$-linear functors $\mod(R)\ra\mod(K)$ where $R$ is a locally bounded $K$-category. The Krull-Gabriel dimension $\KG(R)$ of $R$ is defined as the Krull-Gabriel dimension $\KG(\CF(R))$ of $\CF(R)$. We give in Lemma 4.3 some criterion for $\KG(R)$ to be determined by Krull-Gabriel dimensions of finite convex subcategories of $R$. This lemma is crucial in the proof of Theorem 7.3. 

Section 4 recalls some basic facts about the category $\CF(R)$. These results are particular cases of the very general ones proved in \cite{Au0}, see especially Section 2. Note that in \cite[IV.6]{AsSiSk} the authors show similar facts (in concise and simple manner), but mostly assume that $R$ is a finite dimensional $K$-algebra (equivalently, a finite locally bounded $K$-category). It is easy to see that proofs given in \cite[IV.6]{AsSiSk} are valid in the general setting of locally bounded $K$-categories as well. Therefore we view \cite[IV.6]{AsSiSk} as a good reference for the basic information on the category $\CF(R)$.

Sections 5 and 6 are the core of the paper. Assume that $R$ is a locally support-finite locally bounded $K$-category, $G$ a torsion-free admissible group of $K$-linear automorphisms of $R$ and $F:R\ra A\cong R\slash G$ the Galois covering. Section 5 is devoted to some exact functor $\Phi:\CF(R)\ra\CF(A)$. The definition of $\Phi:\CF(R)\ra\CF(A)$ uses the fact that $F_{\lambda}:\mod(R)\ra\mod(A)$ is a Galois covering of module categories, see Theorem 2.1. In Theorem 5.5 we prove the main properties of $\Phi:\CF(R)\ra\CF(A)$. Section 6 contains the main results of the paper. We prove Theorem 6.3 which is a more detailed version of Theorem 1.2, asserting some further properties of the functor $\Phi:\CF(R)\ra\CF(A)$. In Corollary 6.4 we show that the functor $\Phi$ preserves Krull-Gabriel filtrations, in some concrete sense, provided $\Phi$ is dense. Proofs of Theorem 6.3 and Corollary 6.4 rely on results of Section 5. 

Sections 7 and 8 are devoted to the applications of the main results. Namely, in Section 7 we determine the Krull-Gabriel dimension of locally support-finite repetitive $K$-categories. This is done in Theorem 7.3 by using Lemma 4.3, classification of tame locally support-finite $K$-categories \cite{AsSk4} and classification of cycle-finite algebras with finite Krull-Gabriel dimension \cite{Sk5}. Section 8 applies Theorem 6.3, Theorem 7.3 and the results of \cite{Sk0}, \cite{Sk4} to determine the Krull-Gabriel dimension of standard selfinjective algebras of polynomial growth, see Theorem 8.1. We deduce in Corollary 8.2 that if $A$ is such an algebra of infinite representation type, then $\KG(A)=2$ if and only if $\rad_{A}^{\omega}$ is nilpotent. Finally, we discuss the existence of super-decomposable pure-injective modules over these algebras, see Theorem 8.3.

\section{Galois coverings of locally bounded $K$-categories}

Throughout the paper, $K$ is a fixed algebraically closed field. We denote by $\Mod(K)$ and $\mod(K)$ the categories of all $K$-vector spaces and all finite dimensional $K$-vector spaces, respectively. By an \textit{algebra} we mean a finite dimensional associative basic $K$-algebra with a unit. Assume that $A$ is an algebra. By an \textit{$A$-module} we mean right $A$-module. We denote by $\Mod(A)$ and $\mod(A)$ the categories of all $A$-modules and all finitely generated $A$-modules, respectively. 

Assume that $R$ is a $K$-category in which distinct objects are not isomorphic. Then the objects of $R$ are denoted by lower case letters. Let $\ob(R)$ be the class of all objects of $R$. If $x,y\in\ob(R)$, then $R(x,y)$ denotes the space of all morphisms from $x$ to $y$. Assume that $R$ is a $K$-category. Following \cite{Ga}, \cite{BoGa}, we say that $R$ is \textit{locally bounded} if and only if 

\begin{enumerate}[\rm(1)]
	\item distinct objects of $R$ are not isomorphic,
	\item the algebra $R(x,x)$ is local, for any $x\in\ob(R)$,
	\item $\sum_{y\in\ob(R)}\dim_{K}R(x,y)<\infty$, $\sum_{y\in\ob(R)}\dim_{K}R(y,x)<\infty$, for any $x\in\ob(R)$.
\end{enumerate}

An important example of a locally bounded $K$-category is the \textit{bound quiver} $K$-category $\und{KQ}\slash I$ associated to a bound quiver $(Q,I)$. We briefly recall the definition of a bound quiver $K$-category.

Assume that $Q=(Q_{0},Q_{1})$ is a quiver with the set $Q_{0}$ of vertices and the set $Q_{1}$ of arrows. We say that $Q$ is \textit{finite} if and only if the sets $Q_{0},Q_{1}$ are finite. If $\alpha\in Q_{1}$, then we denote by $s(\alpha)$ its starting vertex and by $t(\alpha)$ its terminal vertex. Assume that $x,y\in Q_{0}$. By a \textit{path} from $x$ to $y$ in $Q$ we mean a sequence $c_{1}...c_{n}$ in $Q_{1}$ such that $s(c_{1})=x$, $t(c_{n})=y$ and $t(c_{i})=s(c_{i+1})$  for $1\leq i<n$. We associate the \textit{stationary path} $e_{x}$ to each vertex $x\in Q_{0}$ and we set $s(e_{x})=t(e_{x})=x$.

Assume that $Q$ is a \textit{locally finite} quiver, that is, the number of arrows in $Q_{1}$ starting or ending in any point is finite. Then the \textit{path $K$-category} $\und{KQ}$ is a $K$-category whose objects are the vertices of $Q$ and the $K$-linear space of morphisms from $x$ to $y$ is generated by all paths from $x$ to $y$. The composition in $\und{KQ}$ is defined by concatenation of paths in $Q$. If $I$ is an admissible ideal in $\und{KQ}$ (see \cite{Pog}), then the pair $(Q,I)$ is called the \textit{bound quiver} and the associated quotient $K$-category $\und{KQ}\slash I$ is the \textit{bound quiver $K$-category}. It is easy to see that this is a locally bounded $K$-category. In fact, every locally bounded $K$-category is isomorphic to some bound quiver $K$-category, see \cite{BoGa}.

Assume that $R$ is a $K$-category (not necessarily locally bounded). A full subcategory $B$ of $R$ is \textit{convex} if and only if for any $n\geq 1$ and objects $x,z_{1},...,z_{n},y$ of $R$ the following condition is satisfied: if $x,y$ are objects of $B$ and the vector spaces of morphisms $R(x,z_{1})$, $R(z_{1},z_{2})$,...,$R(z_{n-1},z_{n})$, $R(z_{n},y)$ are nonzero, then $z_{1},...,z_{n}$ are objects of $B$. Observe that if $R=\und{KQ}/I$ is the path $K$-category of $(Q,I)$ and $B$ is a convex subcategory of $R$, then $B=\und{KQ'}/I'$ for some convex subquiver $Q'$ of the quiver $Q$. Recall that a full subquiver $Q'$ of $Q$ is \textit{convex} if and only if for any path $c_{1}...c_{n}$ from the vertex $x$ to the vertex $y$ in $Q$ such that $x,y\in Q'_{0}$ we have $s(c_{i})\in Q'_{0}$, for $i=1,...,n-1$.

Assume that $R$ is a locally bounded $K$-category. A \textit{right $R$-module} (or an \textit{$R$-module}) is a $K$-linear contravariant functor from $R$ to the category $\Mod(K)$ of $K$-vector spaces. An $R$-module $M:R^{\op}\ra\Mod(K)$ is \textit{finite dimensional} if and only if $\dim M=\sum_{x\in\ob(R)}\dim_{K} M(x)<\infty$. The categories of all $R$-modules and all finite dimensional $R$-modules are denoted by $\Mod(R)$ and $\mod(R)$, respectively. The full subcategories of $\Mod(R)$ and $\mod(R)$ formed by all indecomposable modules are denoted by $\Ind(R)$ and $\ind(R)$, respectively. 

Assume that $R$ is a locally bounded $K$-category and $M,N:R^{\op}\ra\Mod(K)$ are $R$-modules. An \textit{$R$-module homomorphism} $f:M\ra N$ is a natural transformation of functors $(f_{x})_{x\in\ob(R)}$ where $f_{x}:M(x)\ra N(x)$ is a vector space homomorphism, for any $x\in\ob(R)$. The space of all homomorphisms from $M$ to $N$ is denoted by $\Hom_{R}(M,N)$. We usually write ${}_{R}(M,N)$ instead of $\Hom_{R}(M,N)$.

Assume that $R=\und{KQ}\slash I$ is a bound quiver $K$-category. Observe that modules in $\Mod(R)$ can be identified with $K$-linear representations of the bound quiver $(Q,I)$. Similarly, the $R$-module homomorphisms can be identified with morphisms of representations of the bound quiver $(Q,I)$. Furthermore, if $(Q,I)$ is finite, then there is an equivalence of $K$-categories of all $R$-modules and all modules over the \textit{bound quiver $K$-algebra} $KQ\slash I$, see Chapters II and III of \cite{AsSiSk}. This equivalence restricts to an equivalence of categories of finite dimensional modules. Here, the admissible ideal in the $K$-category is denoted by the same letter $I$ as the corresponding admissible ideal in the $K$-algebra. If $(Q,I)$ is finite, we identify the bound quiver $K$-category $R=\und{KQ}\slash I$ with the bound quiver $K$-algebra $KQ\slash I$. Note that any algebra $A$ is isomorphic with some bound quiver $K$-algebra, see Chapter II of \cite{AsSiSk}. If $A\cong KQ\slash I$, then we view $A$-modules as $K$-linear representations of $(Q,I)$.



We recall from \cite{BoGa} the notion of a Galois covering functor, as well as some related concepts. We refer to \cite{Ga} for a general definition of a covering functor and related functors between module categories.

Assume that $R,A$ are locally bounded $K$-categories, $F:R\ra A$ is a $K$-linear functor and $G$ a group of $K$-linear automorphisms of $R$ acting freely on the objects of $R$. This means that $gx=x$ if and only if $g=1$, for any $g\in G$ and $x\in\ob(R)$. Then $F:R\ra A$ is a \textit{Galois covering} if and only if 
\begin{enumerate}[\rm(1)]
	\item the functor $F:R\ra A$ induces isomorpisms $$\bigoplus_{g\in G}R(gx,y)\cong A(F(x),F(y))\cong\bigoplus_{g\in G}R(x,gy)$$ of vector spaces, for any $x,y\in\ob(R)$,
	\item the functor $F:R\ra A$ induces a surjective function $\ob(R)\ra\ob(A)$,
	\item $Fg=F$, for any $g\in G$,
	\item for any $x,y\in\ob(R)$ such that $F(x)=F(y)$ there is $g\in G$ such that $gx=y$. 
\end{enumerate}
We recall that a functor $F:R\ra A$ satisfies the above conditions if and only if $F$ induces an isomorphism $A\cong R\slash G$ where $R\slash G$ is the \textit{orbit category}, see \cite{BoGa}.

Assume that $F:R\ra A\cong R\slash G$ is a Galois covering. The \textit{pull-up} functor $F_{\bullet}:\Mod(A)\ra\Mod(R)$ associated with $F$ is the functor $(-)\circ F^{\op}$. The pull-up functor has the left adjoint $F_{\lambda}:\Mod(R)\ra\Mod(A)$ and the right adjoint $F_{\rho}:\Mod(R)\ra\Mod(A)$ which are called the \textit{push-down} functors. Since the push-down functor $F_{\lambda}$ plays an important role in the paper, we recall its description below. 

Assume that $M:R^{\op}\ra\Mod(K)$ is an $R$-module. We define the $A$-module $F_{\lambda}(M):A^{\op}\ra\Mod(K)$ as follows. Assume that $a\in\ob(A)$ and $a=F(x)$, for some $x\in\ob(R)$. Then $F_{\lambda}(M)(a)=\bigoplus_{g\in G}M(gx)$. Assume that $\alpha\in A(b,a)$ and $a=F(x),b=F(y)$, for some $x,y\in\ob(R)$. Since $F$ induces an isomorphism $\bigoplus_{g\in G}R(gy,x)\cong A(F(y),F(x))$, there are $\alpha_{g}:gy\ra x$, for $g\in G$, such that $\alpha=\sum_{g\in G}F(\alpha_{g})$. Then the homomorphism $F_{\lambda}(M)(\alpha):F_{\lambda}(M)(a)\ra F_{\lambda}(M)(b)$ is defined by homomorphisms $M(g\alpha_{g^{-1}h}):M(gx)\ra M(hy)$, for any $g,h\in G$. \footnote{We use the standard matrix notation for homomorphisms between finite direct sums, see Section 3 for details.} Note that $F_{\rho}(M)(a)=\prod_{g\in G}M(gx)$ and the rest of the definition is similar to the case of $F_{\lambda}$. It follows easily that $F_{\lambda}$ equals $F_{\rho}$ on the category of finite dimensional $A$-modules. We emphasize that this property does not hold in general (but it holds in the case $F:R\ra A$ is a Galois coverings).

Assume that $f:M\ra N$ is an $R$-module homomorphism and $f=(f_{x})_{x\in\ob(R)}$, $f_{x}:M(x)\ra N(x)$. Then $F_{\lambda}(f):F_{\lambda}(M)\ra F_{\lambda}(N)$, $F_{\lambda}(f)=(\hat{f}_{a})_{a\in\ob(A)}$ and $\hat{f}_{a}:F_{\lambda}(M)(a)\ra F_{\lambda}(N)(a)$ is defined by homomorphisms $f_{gx}:M(gx)\ra N(gx)$, for any $g\in G$.  

Observe that if an $R$-module $M$ is finite dimensional, then $F_{\lambda}(M)$ is finite dimensional. Hence the functor $F_{\lambda}$ restricts to a functor $\mod(R)\ra\mod(A)$. This functor is also denoted by $F_{\lambda}$. 

Assume that $R$ is a locally bounded $K$-category, $G$ is a group of $K$-linear automorphisms of $R$ acting freely on the objects of $R$ and $g\in G$. Given $R$-module $M$ we denote by ${}^{g}M$ the module $M\circ g^{-1}$. Given $R$-module homomorphism $f:M\ra N$ we denote by ${}^{g}f$ the $R$-module homomorphism ${}^{g}M\ra {}^{g}N$ such that ${}^{g}f_{x}=f_{g^{-1}x}$, for any $x\in\ob(R)$. It is easy to see that this defines an action of $G$ on $\Mod(R)$. Moreover, the map $f\mapsto {}^{g}f$ defines isomorphism of vector spaces ${}_{R}(M,N)\cong{}_{R}({}^{g}M,{}^{g}N)$.

Assume that $F:R\ra A\cong R\slash G$ is the Galois covering. Recall from \cite{BoGa} that $F_{\lambda}({}^{g}M)\cong F_{\lambda}(M)$, $F_{\lambda}({}^{g}f)=F_{\lambda}(f)$ and $F_{\bullet}(F_{\lambda}(M))\cong\bigoplus_{g\in G}{}^{g}M$, for any $R$-module $M$, $R$-homomorphism $f$ and $g\in G$. Assume that $X,Y\in\mod(R)$. There are only finitely many $g\in G$ such that ${}_{R}(X,{}^{g}Y)\neq 0$, because $G$ acts freely on the objects of $R$. Since the pair $(F_{\lambda},F_{\bullet})$ is an adjoint pair, we get the following isomorphisms of vector spaces $${}_{A}(F_{\lambda}(X),F_{\lambda}(Y))\cong\bigoplus_{g\in G}{}_{R}(X,{}^{g}Y)\cong\bigoplus_{g\in G}{}_{R}({}^{g}X,Y),$$ for any $X,Y\in\mod(R)$. These isomorphisms are natural in both variables $X,Y$ and hence the bifunctors ${}_{A}(F_{\lambda}(-),F_{\lambda}(\cdot)),\bigoplus_{g\in G}{}_{R}(-,{}^{g}(\cdot)),\bigoplus_{g\in G}{}_{R}({}^{g}(-),\cdot)$ are isomorphic. 

Assume that $R$ is a locally bounded $K$-category and $G$ is a group of $K$-linear automorphisms of $R$. Then $G$ is \textit{admissible} if and only if $G$ acts freely on the objects of $R$ and there are only finitely many $G$-orbits. In this case, the orbit category $R\slash G$ is finite and we treat it as an algebra. Assume that $G$ is admissible. We say that $G$ \textit{acts freely on $\ind(R)$} if and only if ${}^{g}M\cong M$ implies that $g=1$, for any $M\in\ind(R)$ and $g\in G$. Assume that $M,L\in\ind(R)$. If $G$ acts freely on $\ind(R)$, then $F_{\lambda}(M)\in\ind(A)$ and $F_{\lambda}(M)\cong F_{\lambda}(L)$ yields $L\cong {}^{g}M$, for some $g\in G$ (see for example \cite{MP}). Moreover, if $G$ is torsion-free, then $G$ acts freely on $\ind(R)$.



Assume that $R$ is a locally bounded $K$-category and $M\in\Mod(R)$. The \textit{support} $\supp (M)$ of $M$ is the full subcategory of $R$ formed by all objects $x$ in $R$ such that $M(x)\neq 0$. The category $R$ is \textit{locally support-finite} if and only if for any $x\in R$ the union of the sets $\supp(M)$, where $M\in\ind(R)$ and $M(x)\neq 0$, is finite. If $R$ is locally support-finite, $G$ is admissible group of $K$-linear automorphisms of $R$ and $F:R\ra A\cong R\slash G$ is the Galois covering, then the results of \cite{DoSk} yield that the push-down functor $F_{\lambda}:\mod(R)\ra\mod(A)$ is dense.

We summarize the above remarks in the following theorem, see \cite[2.5]{DoSk}. This theorem plays a fundamental role in our main results.

\begin{thm} Assume that $R$ is a locally support-finite $K$-category, $G$ an admissible torsion-free group of $K$-linear automorphisms of $R$ and $F:R\ra A$ the Galois covering. Then the functor $F_{\lambda}:\mod(R)\ra\mod(A)$ satisfies the following assertions.
\begin{enumerate}[\rm(1)]
	\item The functor $F_{\lambda}$ induces the following isomorphisms of vector spaces $$\bigoplus_{g\in G}{}_{R}({}^{g}X,Y)\cong{}_{A}(F_{\lambda}(X),F_{\lambda}(Y))\cong\bigoplus_{g\in G}{}_{R}(X,{}^{g}Y),$$ for any $X,Y\in\mod(R)$.
	\item The functor $F_{\lambda}$ is dense, that is, for any $M\in\mod(A)$ there is $X\in\mod(R)$ such that $F_{\lambda}(X)\cong M$.
	\item Assume that $X\in\mod(R)$. Then $F_{\lambda}(X)\cong F_{\lambda}({}^{g}X)$, for any $g\in G$.
	\item Assume that $X,Y\in\ind(R)$. Then $F_{\lambda}(X)\in\ind(A)$ and $F_{\lambda}(X)\cong F_{\lambda}(Y)$ yields $Y\cong {}^{g}X$ for some $g\in G$.
\end{enumerate}
\end{thm}

We recall from \cite{DoSk} (see also \cite{BoGa}) that if the assumptions of Theorem 2.1 hold, then the functor $F_{\lambda}:\mod(R)\ra\mod(A)$ preserves right and left minimal almost split homomorphisms and Auslander-Reiten sequences (see \cite[IV]{AsSiSk} for basic notions of Auslander-Reiten theory). Moreover, the functor $F_{\lambda}$ induces the isomorphism $\ind(R)\slash G\cong\ind(A)$.

We emphasize similarities between the properties of $F_{\lambda}$ from the above theorem and the definition of a Galois covering of $K$-categories. Following \cite{DoSk}, we say that the functor $F_{\lambda}:\mod(R)\ra\mod(A)$ is a \textit{Galois covering of module categories} (in case $R$ and $G$ are as in Theorem 2.1).

More generally, assume that $\CC,\CD$ are additive categories and $F:\CC\ra\CD$ is an additive functor (see \cite[VIII.2]{McL}). Assume that $G$ is a group acting freely on the isomorphism classes of indecomposable objects of $\CC$. We call $F:\CC\ra\CD$ a \textit{Galois covering of additive categories} if and only if $F$ satisfies analogous assertions as the push-down functor from Theorem 2.1. Note that in this case $\CD$ is isomorphic (as an additive category) to the orbit category $\CC\slash G$.

\section{Krull-Gabriel dimension of abelian categories}


Throughout the paper we assume that any abelian category $\CC$ is \textit{skeletally small}. This means that the class of all isomorphism classes of objects of $\CC$ is a set.

Assume that $\CC$ is an abelian category and $X,Y$ are objects of $\CC$. We denote by $\Hom_{\CC}(X,Y)$ the abelian group of all homomorphisms from $X$ to $Y$. We usually write ${}_{\CC}(X,Y)$ instead of $\Hom_{\CC}(X,Y)$. If $X'$ is a subobject of $X$, then we write $X'\subseteq X$ and in this case $X\slash X'$ denotes the quotient object. If $f\in{}_{\CC}(X,Y)$, then $\Ker_{\CC}(f)$, $\Im_{\CC}(f)$, $\Coker_{\CC}(f)$ and $\Coim_{\CC}(f)$ denote \textit{kernel}, \textit{image}, \textit{cokernel} and \textit{coimage} of $f$ in $\CC$, respectively. It is well-known that these objects are some special subobjects or quotients of $X$ and $Y$. We refer to \cite{Po} for precise definitions of these notions in abelian categories. 

Assume that $\CC$ is an abelian category and $X_{1},...,X_{n},Y_{1},...,Y_{m}$ are objects of $\CC$. We use the standard \textit{matrix notation} for $f\in{}_{\CC}(\bigoplus_{i=1}^{n}X_{i},\bigoplus_{j=1}^{m}Y_{j})$. This means that we write $f=[f_{ji}]_{i=1,...,n}^{j=1,...,m}$ where $p_{j}:\bigoplus_{k=1}^{n}Y_{k}\ra Y_{j}$ is the split epimorphism, $u_{i}:X_{i}\ra\bigoplus_{k=1}^{m}X_{k}$ is the split monomorphism and $f_{ji}=p_{j}fu_{i}$, for any $i=1,...,n$, $j=1,...,m$. Moreover, we say in this case that $f$ \textit{is defined by homomorphisms $f_{ji}$}.

Assume that $\CC,\CD$ are abelian categories. A functor $F:\CC\ra\CD$ is \textit{exact} if and only if for any exact sequence $0\ra X\stackrel{f}{\rightarrow}Z\stackrel{g}{\rightarrow}Y\ra 0$ in $\CC$ the sequence $0\ra F(X)\stackrel{F(f)}{\longrightarrow}F(Z)\stackrel{F(g)}{\longrightarrow}F(Y)\ra 0$ is exact in $\CD$. Equivalently, the functor $F$ preserves finite limits and finite colimits, see \cite[VIII]{McL}. Note that if $F:\CC\ra\CD$ is exact, then $X'\subseteq X$ yields $F(X')\subseteq F(X)$ and $F(X\slash X')=F(X)\slash F(X')$, for any objects $X',X$ of $\CC$. Moreover, $F$ preserves kernels, images, cokernels and coimages of homomorphisms in $\CC$.

Assume that $\CC$ is an abelian category. A full subcategory $\CS$ of $\CC$ is a \textit{Serre subcategory} if and only if for any exact sequence $0\ra X\ra Z\ra Y\ra 0$ in $\CS$ we have $Z\in\CS$ if and only if $X,Y\in\CS$. Therefore $\CS$ is closed under subobjects, quotients and extensions. It is easy to see that a Serre subcategory is an abelian subcategory. 

Basic examples of Serre subcategories are provided by kernels of exact functors. Indeed, if $F:\CC\ra\CD$ is an exact functor, then the \textit{kernel} $\Ker(F)$ of $F$ is the full subcategory of $\CC$ formed by all objects $X\in\CC$ such that $F(X)=0$. It is easy to see that $\Ker(F)$ is a Serre subcategory of $\CC$.

Assume that $\CS$ is a Serre subcategory of $\CC$. We recall the definition of the \textit{quotient category} $\CC\slash\CS$. The class of objects of $\CC\slash\CS$ coincides with the class of objects of $\CC$. Homomorphisms in $\CC\slash\CS$ are defined in the following way. Assume that $X,Y$ are objects of $\CC$ and denote by $\CS_{X,Y}$ the set of all abelian groups ${}_{\CC}(X',Y\slash Y')$ where $X'\subseteq X$, $Y'\subseteq Y$ and $X\slash X',Y'\in\CS$. Then $\CS_{X,Y}$ is a directed system of abelian groups, see \cite{Po}, and ${}_{\CC\slash\CS}(X,Y)$ is defined as the direct limit $\varinjlim{}_{\CC}(X',Y\slash Y')$ where ${}_{\CC}(X',Y\slash Y')\in\CS_{X,Y}$. The quotient category $\CC\slash\CS$ is abelian and there exists an exact functor $q_{\CS}:\CC\ra\CC\slash\CS$ such that $q_{\CS}(X)=X$ for any object $X$ of $\CC$. We refer to \cite{Po} for the precise definition of $q_{\CS}$ on homomorphisms. Nevertheless, we recall from Section 4.3 of \cite{Po} that if $f$ is a homomorphism in $\CC$, then $q_{\CS}(f)=0$ if and only if $\Im_{\CC}(f)\in\CS$, $q_{\CS}(f)$ is a monomorphism if and only if $\Ker_{\CC}(f)\in\CS$ and $q_{\CS}(f)$ is an epimorphism if and only if $\Coker_{\CC}(f)\in\CS$. The functor $q_{\CS}:\CC\ra\CC\slash\CS$ is called the \textit{quotient functor}. 



Assume that $\CC$ is an abelian category. An object $S$ of $\CC$ is \textit{simple} in $\CC$ if and only if $S\neq 0$ and any subobject of $S$ in $\CC$ is either $S$ or $0$. An object $T$ of $\CC$ has \textit{finite length} in $\CC$ if and only if there exists a chain of subobjects $0=T_{0}\subseteq T_{1}\subseteq ...\subseteq T_{n}=T$ of $T$ such that $T_{i+1}\slash T_{i}$ is simple in $\CC$, for any $i=0,...,n-1$. We call such a chain a \textit{composition series} of $T$. It is well-known that if $T$ has finite length, then the number $n$ is unique. We call the number $n$ the \textit{length} of $T$ and denote it by $l_{\CC}(T)$. Recall that if $0\ra X\ra Z\ra Y\ra 0$ is a short exact sequence in $\CC$, then $l_{\CC}(Z)=l_{\CC}(X)+l_{\CC}(Y)$ and hence $Z$ has finite length if and only if $X,Y$ have finite length.

Assume that $\CC$ is an abelian category and $\CS$ is a Serre subcategory of $\CC$. An object $S$ of $\CC$ is \textit{$\CS$-simple} if and only if $S\notin\CS$ and $S'\in\CS$ or $S\slash S'\in\CS$, for any subobject $S'$ of $S$ in $\CC$. The following fact is applied in Section 6.

\begin{lm} Assume that $\CC$ is an abelian category and $\CS$ a Serre subcategory of $\CC$. The following assertions hold.
\begin{enumerate}[\rm(1)]
	\item An object $S$ of $\CC$ is simple in $\CC\slash\CS$ if and only if $S$ is $\CS$-simple.
	\item An object $T$ of $\CC$ has finite length in $\CC\slash\CS$ if and only if there exists a chain of subobjects $0=T_{0}\subseteq T_{1}\subseteq ...\subseteq T_{n}=T$ of $T$ in $\CC$ such that $T_{i+1}\slash T_{i}$ is $\CS$-simple, for any $i=0,...,n-1$.
\end{enumerate}
\end{lm}

{\bf Proof.} (1) Assume that $S$ is an object of $\CC$ and $S$ is simple in $\CC\slash\CS$. We show that $S$ is $\CS$-simple. Assume that $S\in\CS$ and let $1_{S}:S\ra S$ be the identity homomorphism of $S$ in $\CC$. Then $\Im_{\CC}(1_{S})=S\in\CS$ and thus $q_{\CS}(1_{S})=0$. Since $q_{\CS}(1_{S})$ is also the identity of $S$ in $\CC\slash\CS$, we get $S\cong 0$ in $\CC\slash\CS$, contradiction. This shows that $S\notin\CS$. 

Assume that $S'$ is a subobject of $S$ in $\CC$. Then there exists a short exact sequence $0\ra S'\stackrel{f}{\ra} S\stackrel{g}{\ra} S\slash S'\ra 0$ in $\CC$ which yields a short exact sequence $$0\ra q_{\CS}(S')\stackrel{q_{\CS}(f)}{\longrightarrow} q_{\CS}(S)\stackrel{q_{\CS}(g)}{\longrightarrow} q_{\CS}(S)\slash q_{\CS}(S')\ra 0$$ in $\CC\slash\CS$. Since $q_{\CS}(S)=S$ is simple in $\CC\slash\CS$, we get that $q_{\CS}(f)$ or $q_{\CS}(g)$ is an isomorphism in $\CC\slash\CS$. Thus $q_{\CS}(f)$ is an epimorphism or $q_{\CS}(g)$ is a monomorphism and hence $\Coker_{\CC}(f)=S\slash S'\in\CS$ or $\Ker_{\CC}(g)=S'\in\CS$. Similar arguments show the other implication. 

(2) The assertion follows from (1) and the fact that the quotient functor is an identity on objects. \epv 

Assume that $\CC$ is an abelian category. Following \cite{Ge2} (see also \cite{Ga0}), we define the \textit{Krull-Gabriel filtration} $(\CC_{\alpha})_{\alpha}$ of $\CC$ recursively as follows: 
\begin{enumerate}[\rm(1)]
	\item $\CC_{-1}=0$,
	\item $\CC_{\alpha+1}$ is the Serre subcategory of $\CC$ formed by all functors having finite length in the quotient category $\CC\slash\CC_{\alpha}$, for any ordinal number $\alpha$,
	\item $\CC_{\beta}=\bigcup_{\alpha<\beta}\CC_{\alpha}$, for any limit ordinal $\beta$.
\end{enumerate} 

Observe that $\CC_{0}$ is the Serre subcategory of $\CC$ formed by all objects of $\CC$ of finite length. It is easy to see that $\CC_{\alpha}$ is a Serre subcategory of $\CC_{\alpha+1}$, for any ordinal $\alpha\geq -1$.

Assume that $(\CC_{\alpha})_{\alpha}$ is the Krull-Gabriel filtration of $\CC$. We define the \textit{Krull-Gabriel dimension} $\KG(\CC)$ of $\CC$ as the smallest ordinal number $\alpha$ such that $\CC_{\alpha}=\CC$, if such a number exists. Otherwise, we set $\KG(\CC)=\infty$ and take a technical assumption that $\alpha<\infty$, for any ordinal number $\alpha$. If $\KG(\CC)=n\in\NN$, then the Krull-Gabriel dimension of $\CC$ is \textit{finite}. If $\KG(\CC)=\infty$, then the Krull-Gabriel dimension of $\CC$ is \textit{undefined}.

The following lemma gives a sufficient condition for a functor to preserve the Krull-Gabriel dimension. We apply the lemma in Sections 4 and 6.

\begin{lm} Assume that $\CC,\CD$ are abelian categories, $(\CC_{\alpha})_{\alpha},(\CD_{\alpha})_{\alpha}$ are Krull-Gabriel filtrations of $\CC,\CD$, respectively, and $F:\CC\ra\CD$ is an exact functor.
\end{lm}
\begin{enumerate}[\rm(1)]
	\item Assume that for any $U\in\CD$ there exists an epimorphism $\epsilon:F(T)\ra U$, for some $T\in\CC$. Assume that $T\in\CC_{\alpha}$ implies $F(T)\in\CD_{\alpha}$, for any $T\in\CC$ and any ordinal $\alpha$. Then $\KG(\CD)\leq\KG(\CC)$.
	\item Assume that $F(T)\in\CD_{\alpha}$ implies $T\in\CC_{\alpha}$, for any $T\in\CC$ and any ordinal $\alpha$. Then $\KG(\CC)\leq\KG(\CD)$.
	\item Assume that for any $U\in\CD$ there exists an epimorphism $\epsilon:F(T)\ra U$, for some $T\in\CC$. Assume that $T\in\CC_{\alpha}$ if and only if $F(T)\in\CD_{\alpha}$, for any $T\in\CC$ and any ordinal $\alpha$. Then $\KG(\CC)=\KG(\CD)$.
\end{enumerate}

{\bf Proof.} (1) Assume that $\alpha$ is an ordinal number and $\alpha+1\leq\KG(\CD)$. Then $\CD_{\alpha}\neq\CD$ and hence there is $U\in\CD$ such that $U\notin\CD_{\alpha}$. Since there is an epimorphism $\epsilon:F(T)\ra U$, for some $T\in\CC$, we get $F(T)\notin\CD_{\alpha}$ and thus $T\notin\CC_{\alpha}$. This implies that $\CC_{\alpha}\neq\CC$, hence $\alpha+1\leq\KG(\CC)$. It follows that $\beta\leq\KG(\CD)$ yields $\beta\leq\KG(\CC)$, for any ordinal number $\beta$ (both successor ordinal and limit ordinal). Consequently we get $\KG(\CD)\leq\KG(\CC)$.

(2) Assume that $\alpha$ is an ordinal number and $\alpha+1\leq\KG(\CC)$. Then $\CC_{\alpha}\neq\CC$ and hence there is $T\in\CC$ such that $T\notin\CC_{\alpha}$. This yields $F(T)\notin\CD_{\alpha}$, so $\CD_{\alpha}\neq\CD$ and $\alpha+1\leq\KG(\CD)$. It follows that $\beta\leq\KG(\CC)$ yields $\beta\leq\KG(\CD)$, for any ordinal number $\beta$ (both successor ordinal and limit ordinal). Consequently we get $\KG(\CC)\leq\KG(\CD)$.

(3) The assertion follows directly from (1) and (2). \epv

Assume that $\CC,\CD$ are abelian categories, $(\CC_{\alpha})_{\alpha},(\CD_{\alpha})_{\alpha}$ are Krull-Gabriel filtrations of $\CC,\CD$, respectively, and $F:\CC\ra\CD$ is an exact functor. It follows from Appendix B of \cite{Kr} that if $F:\CC\ra\CD$ is full and dense, then $T\in\CC_{\alpha}$ implies $F(T)\in\CD_{\alpha}$, for any $T\in\CC$ and any ordinal $\alpha$. Moreover, if $F:\CC\ra\CD$ is full and faithul, then $F(T)\in\CD_{\alpha}$ implies $T\in\CC_{\alpha}$, for any $T\in\CC$ and any ordinal $\alpha$. 

Assume that $\CS$ is a Serre subcategory of $\CC$. The quotient functor $q_{\CS}:\CC\ra\CC\slash\CS$ is full, dense and exact, so Lemma 3.2 (1) yields $\KG(\CC\slash\CS)\leq\KG(\CC)$. The inclusion functor $\CS\hookrightarrow\CC$ is full, faithful and exact, so Lemma 3.2 (2) yields $\KG(\CS)\leq\KG(\CC)$. Moreover, it follows from Section 1 of \cite{Ge2} that $\CC$ has finite Krull-Gabriel dimension if and only if $\CS$ and $\CC\slash\CS$ have finite Krull-Gabriel dimensions and there are inequalities $$\KG(\CS),\KG(\CC\slash\CS)\leq\KG(\CC)\leq\KG(\CS)+\KG(\CC\slash\CS)+1.$$





\section{Finitely presented functors and Krull-Gabriel dimension}

Assume that $R$ is a locally bounded $K$-category. We denote by $\CG(R)$ the category of all contravariant $K$-linear functors from $\mod(R)$ to the category $\mod(K)$ of finite dimensional $K$-vector spaces.


Assume that $F,G,H\in\CG(R)$ and let $F\stackrel{u}{\rightarrow}G\stackrel{v}{\rightarrow}H$ be a sequence of homomorphisms of functors. Recall that this sequence is exact in $\CG(R)$ if and only if it induces an exact sequence $F(X)\stackrel{u_{X}}{\rightarrow}G(X)\stackrel{v_{X}}{\rightarrow}H(X)$ of vector spaces, for any $X\in\mod(R)$. If $F=0$, then $v:G\ra H$ is a monomorphism of functors, that is, $v_{X}:G(X)\ra H(X)$ is a monomorphism of vector spaces, for any $X\in\mod(R)$. In this case, $G$ is a subobject of $H$. If $H=0$, then $u:F\ra G$ is an epimorphism of functors, that is, $v_{X}:F(X)\ra G(X)$ is an epimorphism of vector spaces, for any $X\in\mod(R)$. In this case, $G$ is a quotient of $F$.

Assume that $M\in\mod(R)$. Then a \textit{contravariant hom-functor} is the functor $H_{M}:\mod(R)\ra\mod(K)$ such that $H_{M}(X)={}_{R}(X,M)$, for any $X\in\mod(R)$, and if $f\in{}_{R}(X,Y)$, then $H_{M}(f):{}_{R}(Y,M)\ra{}_{R}(X,M)$ where $H_{M}(f)(g)=gf$, for any $g\in{}_{R}(Y,M)$. The functor $H_{M}:\mod(R)\ra\mod(K)$ is denoted by ${}_{R}(-,M)$. 

Assume $f\in{}_{R}(M,N)$. The homomorphism $f$ induces a homomorphism of functors ${}_{R}(-,f):{}_{R}(-,M)\ra{}_{R}(-,N)$ such that ${}_{R}(X,f):{}_{R}(X,M)\ra{}_{R}(X,N)$ is defined by ${}_{R}(X,f)(g)=fg$, for any $g\in{}_{R}(X,M)$. The Yoneda lemma implies that the function $f\mapsto{}_{R}(-,f)$ defines an isomorphism ${}_{R}(M,N)\ra{}_{\CG(R)}({}_{R}(-,M),{}_{R}(-,N))$ of vector spaces. Moreover, this yields $M\cong N$ if and only if ${}_{R}(-,M)\cong{}_{R}(-,N)$. 

Assume that $F\in\CG(R)$. The functor $F$ is \textit{finitely generated} if and only if there exists an epimorphism of functors ${}_{R}(-,N)\ra F$, for some $N\in\mod(R)$. The functor $F$ is \textit{finitely presented} if and only if there exists an exact sequence of functors ${}_{R}(-,M)\stackrel{{}_{R}(-,f)}{\longrightarrow}{}_{R}(-,N)\ra F\ra 0$, for some $M,N\in\mod(R)$ and $R$-module homomorphism $f:M\ra N$. This means that $F\cong\Coker{}_{R}(-,f)$ which yields $F(X)$ is isomorphic to the cokernel of the map ${}_{R}(X,f):{}_{R}(X,M)\ra{}_{R}(X,N)$.

We denote by $\CF(R)$ the full subcategory of $\CG(R)$ formed by finitely presented functors. Obviously ${}_{R}(-,M)\in\CF(R)$ for any $M\in\mod(R)$. Moreover, the functor ${}_{R}(-,M)$ is a projective object of the category $\CF(R)$ and any projective object of $\CF(R)$ is a hom-functor, see \cite{Au0}.

Assume that $N\in\mod(R)$ is indecomposable. The functor $\rad_{R}(-,N)$, sending an $R$-module $X$ to the Jacobson radical $\rad_{R}(X,N)$ of ${}_{R}(X,N)$, is the unique maximal subfunctor of the functor ${}_{R}(-,N)$. This implies that the quotient functor $S^{N}={}_{R}(-,N)\slash\rad_{R}(-,N)$ is simple. Moreover, every simple object of the category $\CF(R)$ is isomorphic with a functor $S^{N}$, for some indecomposable $R$-module $N$. Observe that $S^{N}(N)=\End_{R}(N)\slash\rad(\End_{R}(N))$ is a one dimensional vector space and hence there is a unique (up to a scalar multiple) epimorphism of functors $\pi^{N}:{}_{R}(-,N)\ra S^{N}$. Moreover, an $R$-module homomorphism $g:M\ra N$ is a right minimal almost split if and only if the induced sequence of functors ${}_{R}(-,M)\stackrel{{}_{R}(-,g)}{\longrightarrow}{}_{R}(-,N)\stackrel{\pi^{N}}{\ra}S^{N}\ra 0$ is a minimal projective presentation of $S^{N}$ in $\CF(R)$, see Section 2 of \cite{Au0}.

Assume that $R$ is a locally bounded $K$-category. We define the \textit{Krull-Gabriel dimension of $R$} as the Krull-Gabriel dimension $\KG(\CF(R))$ of the category $\CF(R)$ and denote it by $\KG(R)$. 


Assume that $\CX$ is a full subcategory of $\mod(R)$. If $M\in\CX$ is an $R$-module, then ${}_{\CX}(-,M):\CX\ra\mod(K)$ is the functor ${}_{R}(-,M)$ restricted to the category $\CX$. If $M,N\in\CX$ are $R$-modules, then any $R$-homomorphism $f:M\ra N$ induces a homomorphism of functors ${}_{\CX}(-,f):{}_{\CX}(-,M)\ra{}_{\CX}(-,N)$. A contravariant functor $F:\CX^{\op}\ra\mod(K)$ is \textit{finitely generated} if and only if there exists an exact sequence of functors ${}_{\CX}(-,N)\ra F\ra 0$, for some $N\in\CX$. A functor $F:\CX^{\op}\ra\mod(K)$ is \textit{finitely presented} if and only if there exists an exact sequence of functors ${}_{\CX}(-,M)\stackrel{{}_{\CX}(-,f)}{\longrightarrow}{}_{\CX}(-,N)\ra F\ra 0$, for some $M,N\in\CX$ and $R$-module homomorphism $f:M\ra N$. We denote by $\CF(\CX)$ the category of all finitely presented functors $\CX^{\op}\ra\mod(K)$. 

An important role in the study of the Krull-Gabriel dimension $\KG(R)$ of $R$ is played by contravariantly finite subcategories of the category $\mod(R)$, see \cite{AuRe}. Assume that $\CX$ is a full subcategory of $\mod(R)$ which is closed under isomorphisms and direct summands. Then $\CX$ is \textit{contravariantly finite} if and only if any module $M\in\mod(R)$ has a \textit{right $\CX$-approximation}. This means that there exists a module $X_{M}\in\CX$ and $R$-homomorphism $\alpha:X_{M}\ra M$ such that for any module $X\in\CX$ and any $R$-homomorphism $\beta:X\ra M$ there is $\gamma:X\ra X_{M}$ such that $\alpha\gamma=\beta$. This condition is equivalent with the existence of an exact sequence of functors ${}_{\CX}(-,X_{M})\ra {}_{\CX}(-,M)\ra 0$. Thus $\CX$ is contravariantly finite if and only if the functor ${}_{\CX}(-,M)$ is finitely generated, for any $M\in\mod(R)$. If $\CX$ is contravariantly finite, then the category $\CF(\CX)$ is abelian, see \cite{Ge2}.

\begin{prop} Assume that $R$ is a locally bounded $K$-category and $B$ is a convex subcategory of $R$. Then the category $\mod(B)$ is a contravariantly finite subcategory of $\mod(R)$.
\end{prop}

{\bf Proof.} Since $B$ is a convex subcategory of $R$, there exists a functor of extension by zeros $\mathcal{E}:\mod(B)\ra\mod(R)$. This functor is full, faithful and exact, so $\mod(B)$ is a full subcategory of $\mod(R)$. We identify $\mod(B)$ with its image via $\mathcal{E}$.

Assume that $M\in\mod(R)$. We show that $M$ has a right $\mod(B)$-approximation. Indeed, let $M'$ be the $R$-module generated by all modules of the form $\Im(f)\subseteq M$ where $f:X\ra M$ and $X\in\mod(B)$. Then $M'\in\mod(B)$, because $\Im(f)$ is an $B$-module. Moreover, we have $\Im(\beta)\subseteq M'$, for any $X\in\mod(B)$ and $\beta:X\ra M$. This yields the inclusion $M'\hookrightarrow M$ is the right $\mod(B)$-approximation of $M$. \epv 

It is well-known that $\CX$ is a contravariantly finite subcategory of $\mod(R)$ if and only if for any functor $F\in\CF(R)$ the restriction $F_{|\CX}:\CX\ra\mod(K)$ is a finitely presented functor, see for example \cite{Au}. In this case we define the \textit{restriction functor} $r_{\CX}:\CF(R)\ra\CF(\CX)$ such that $r_{\CX}(F)=F_{|\CX}$, for any $F\in\CF(R)$. The restriction functor $r_{\CX}:\CF(R)\ra\CF(\CX)$ is exact, full and dense, see Section 2 of \cite{Ge2}. 

We say that $S\in\CF(R)$ is a \textit{proper subfunctor} of $T\in\CF(R)$ if and only if $S\subseteq T$, $S\neq 0$ and $S$ is not isomorphic with $T$. Assume that $B$ is a convex subcategory of $R$. Proposition 4.1 yields $\CX=\mod(B)$ is a contravariantly finite subcategory of $\mod(R)$. In this case the restriction functor $r_{\CX}:\CF(R)\ra\CF(\CX)$ is denoted by $r_{B}:\CF(R)\ra\CF(B)$. Moreover, we write $\supp(T)\subseteq B$ for any $T\in\CF(R)$ such that $T(X)=0$ for $X\in\ind(R)$ with $\supp(X)\not\subseteq B$. Note that if $U\subseteq T$ and $\supp(T)\subseteq B$, then $\supp(U)\subseteq B$ since $T(X)=0$ yields $U(X)=0$.

\begin{lm} Assume that $R$ is a locally bounded $K$-category and $B$ is a convex subcategory of $R$. Then $r_{B}(T)\in\CF(B)_{\alpha}$ implies $T\in\CF(R)_{\alpha}$ for any $T\in\CF(R)$ such that $\supp(T)\subseteq B$, for any ordinal number $\alpha$.
\end{lm}

{\bf Proof.} We proceed by transfinite induction with respect to an ordinal $\alpha$. Assume that $\alpha=0$ and take $T\in\CF(R)$ such that $\supp(T)\subseteq B$. Assume that $S$ is a proper subfunctor of $T$. We show that $r_{B}(S)$ is a proper subfunctor of $r_{B}(T)$. The functor $r_{B}$ is exact, so $r_{B}(S)\subseteq r_{B}(T)$. Since $S\neq 0$, there is $X\in\ind(R)$ such that $S(X)\neq 0$, hence $T(X)\neq 0$ (because $S(X)\subseteq T(X)$) and thus $\supp(X)\subseteq B$. This yields $r_{B}(S)(X)=S(X)\neq 0$, so $r_{B}(S)\neq 0$. Since $S\subseteq T$ is not an isomorphism, there is $Y\in\ind(R)$ such that $\dim_{K}S(Y)<\dim_{K}T(Y)$. This implies $T(Y)\neq 0$ and thus $\supp(Y)\subseteq B$. Hence we get $r_{B}(S)(Y)=S(Y)$ and $r_{B}(T)(Y)=T(Y)$, so $r_{B}(S)\subseteq r_{B}(T)$ is not an isomorphism. Summing up, $r_{B}(S)$ is a proper subfunctor of $r_{B}(T)$. These arguments imply that if $T$ has arbitrarily long chains of proper subfunctors, then so has $r_{B}(T)$. Consequently, $r_{B}(T)\in\CF(B)_{0}$ yields $T\in\CF(R)_{0}$. 

Assume that $\alpha$ is an ordinal and $r_{B}(T)\in\CF(B)_{\alpha}$ implies $T\in\CF(R)_{\alpha}$ for any $T\in\CF(R)$ such that $\supp(T)\subseteq B$. We show that the same holds for $\alpha+1$. Assume that $\supp(U)\subseteq B$ and $U\in\CF(R)$ is not $\CF(R)_{\alpha}$-simple (i.e. simple in the quotient category $\CF(R)\slash\CF(R)_{\alpha}$, see Lemma 3.1). We show that $r_{B}(U)$ is not $\CF(B)_{\alpha}$-simple. There exists $V\subseteq U$ such that $V\notin\CF(R)_{\alpha}$ and $U\slash V\notin\CF(R)_{\alpha}$. Then $\supp(V),\supp(U\slash V)\subseteq B$, so $r_{B}(V),r_{B}(U\slash V)\notin\CF(B)_{\alpha}$, by the assumption. Hence $r_{B}(U)$ is not $\CF(B)_{\alpha}$-simple, because $r_{B}(V)\subseteq r_{B}(U)$ and $r_{B}(U\slash V)=r_{B}(U)\slash r_{B}(V)$. These arguments imply that if $T$ has no finite length in $\CF(R)\slash\CF(R)_{\alpha}$, then neither has $r_{B}(T)$ in $\CF(B)\slash\CF(B)_{\alpha}$. Consequently, $r_{B}(T)\in\CF(B)_{\alpha+1}$ yields $T\in\CF(R)_{\alpha+1}$. This gives the inductive step for a successor ordinal. Since $\CF(X)_{\beta}=\bigcup_{\alpha<\beta}\CF(X)_{\alpha}$ where $X$ is $R$ or $B$, it is easy to see that the step for a limit ordinal also follows. Hence the assertion follows by transfinite induction. \epv


Assume that $(B_{n})_{n\in\NN}$ is a sequence of locally bounded $K$-categories. If $\KG(B_{n})$ is defined, for any $n\in\NN$, then we set $\sup(\KG(B_{n}))_{n\in\NN}=\bigcup_{n\in\NN}\KG(B_{n})$. Otherwise, we set $\sup(\KG(B_{n}))_{n\in\NN}=\infty$. In Section 7 of the paper we apply the following criterion which is some version of \cite[Corollary 1.5]{Ge2}. 

\begin{lm} Assume that $R$ is a locally support-finite $K$-category and $B_{1}\subseteq B_{2}\subseteq\hdots\subseteq B_{n}\subseteq\hdots\subseteq R$ an ascending chain of finite convex subcategories of $R$ such that $\bigcup_{n\in\NN}B_{n}=R$. Then $\KG(R)=\sup(\KG(B_{n}))_{n\in\NN}$. 
\end{lm}

{\bf Proof.} We introduce some notation. Assume that $x$ is an object of $R$. Then $\CA(x)$ is the full subcategory of $R$ formed by all objects $y$ such that $y\in\supp(N)$, for some indecomposable $R$-module $N$ such that $x\in\supp(N)$. Note that $\CA(x)$ is finite, because $R$ is locally support-finite. If $C$ is some set of objects of $R$, then $\CA(C)$ is the full subcategory of $R$ formed by all objects of $\bigcup_{x\in C}\CA(x)$.

We show that for any $T\in\CF(R)$ there exists $n\in\NN$ such that $\supp(T)\subseteq B_{n}$. Indeed, assume that $T={}_{R}(-,M)$ for some $M\in\mod(R)$. It is easy to see that $T(X)={}_{R}(X,M)\neq 0$ if and only if the support of $X$ is contained in $\CA(\supp(M))$. Therefore $\supp(T)\subseteq\CA(\supp(M))\subseteq B_{n}$, for some $n\in\NN$, because $\CA(\supp(M))$ is finite. Since any functor $T\in\CF(R)$ is a quotient of some hom-functor, the assertion follows.

We show that $\KG(R)=\sup(\KG(B_{n}))_{n\in\NN}$. Since $\CF(B_{n})$ is a Serre subcategory of $\CF(R)$, we get $\KG(B_{n})\leq\KG(R)$, for any $n\in\NN$ (this also follows from the existence of $r_{B_{n}}:\CF(R)\ra\CF(B_{n})$). This yields $\sup(\KG(B_{n}))_{n\in\NN}\leq\KG(R)$. In particular, if $\KG(B_{n})=\infty$, for some $n\in\NN$, then $\KG(R)=\infty$ and so $\KG(R)=\sup(\KG(B_{n}))_{n\in\NN}$ in this case. Hence assume that $\KG(B_{n})$ is defined, for any $n\in\NN$. Assume that $T\in\CF(R)$ and $\supp(T)\subseteq B_{n}$. Set $\alpha=\KG(B_{n})$. Since $r_{B_{n}}(T)\in\CF(B_{n})_{\alpha}$, we get $T\in\CF(R)_{\alpha}$ by Lemma 4.2. These arguments imply that $\KG(R)\leq\sup(\KG(B_{n}))_{n\in\NN}$. Consequently, we obtain $\KG(R)=\sup(\KG(B_{n}))_{n\in\NN}$. \epv

The thesis of the above lemma does not hold without the assumption that $R$ is locally support-finite. Indeed, assume that $R$ is the universal Galois covering of the Kronecker quiver, that is, $R=\underline{KQ}$ where $Q$ is the infinite quiver of the form $$\xymatrix{&\ar[dl]\bullet\ar[dr]&&\ar[dl]\bullet\ar[dr]\\\hdots\bullet&&\bullet&&\bullet\hdots}$$ Observe that every indecomposable $R$-module is linear, that is, a module with $K$ in the vertices and identities on the arrows. It follows easily that $R$ is not locally support-finite. Let $S$ be any simple injective $R$-module and consider the hom-functor $T={}_{R}(-,S)$. It is easy to see that $T(M)\neq 0$ for infinitely many pairwise nonisomorphic indecomposable $R$-modules $M$. This yields $T$ is not of finite length, so $T\notin\CF(R)_{0}$ and thus $\KG(R)\neq 0$ (in fact, it was pointed out by G. Bobi\'nski that using some arguments from \cite{BobKr} one can show that $\KG(R)=2$). On the other hand, every finite convex subcategory of $R$ is of type $\mathbb{A}_{n}$, so its Krull-Gabriel dimension is zero by Auslander's result. Hence the assumption in Lemma 4.3 (that $R$ is locally support-finite) is the crucial one.  

Assume that $R$ is a locally bounded $K$-category. We say that $R$ has \textit{subcategory-determined} Krull-Gabriel dimension if and only if there exists an ascending chain $B_{1}\subseteq B_{2}\subseteq\hdots\subseteq B_{n}\subseteq\hdots\subseteq R$ of finite convex subcategories of $R$ such that $\bigcup_{n\in\NN}B_{n}=R$ and $\KG(R)=\sup(\KG(B_{n}))_{n\in\NN}$. The class of locally bounded $K$-categories with subcategory-determined Krull-Gabriel dimension seems to be important in the study of Krull-Gabriel dimension, just as Lemma 4.3 plays an important role in Section 7.  

\section{The functor $\Phi:\CF(R)\ra\CF(A)$}

Throughout we assume that $R$ is a locally support-finite $K$-category, $G$ a torsion-free admissible group of $K$-linear automorphisms of $R$ and $F:R\ra A\cong R\slash G$ the Galois covering. Recall that in this case the push-down functor $F_{\lambda}:\mod(R)\ra\mod(A)$ is a Galois covering of module categories, see Theorem 2.1. All properties of $F_{\lambda}$ stated in Theorem 2.1, especially the fact that it is dense, are used freely in the section. Moreover, we assume that the pull-up functor $F_{\bullet}$ is restricted to the category $\mod(A)$, that is, we consider $F_{\bullet}:\mod(A)\ra\Mod(R)$.

The section is devoted to define some covariant exact functor $\Phi:\CF(R)\ra\CF(A)$ and show its main properties. This functor is used in Section 6 to prove our main results. We start with two preparatory facts.

\begin{lm} Assume that $T\in\CF(R)$, $T\neq 0$ and $X\in\mod(R)$. Then $T({}^{g}X)\neq 0$ only for finite number of $g\in G$.
\end{lm}

{\bf Proof.} Assume that $T=\Coker{}_{R}(-,f)$ for some $R$-homomorphism $f:M\ra N$. The vector space $T({}^{g}X)$ is a quotient of ${}_{R}({}^{g}X,N)$, for any $g\in G$. Since $G$ acts freely on the objects of $R$ and supports of $X,N$ are finite, there is only finite number of $g\in G$ such that ${}_{R}({}^{g}X,N)\neq 0$. This shows the claim. \epv

The following result is a more concrete description of the isomorphism of bifunctors $\bigoplus_{g\in G}{}_{R}({}^{g}(-),\cdot)\cong{}_{A}(F_{\lambda}(-),F_{\lambda}(\cdot))$, see Section 2 and \cite{Ga}, \cite{BoGa} for the details. In this description we identify $F_{\lambda}({}^{g}X)$ with $F_{\lambda}(X)$, for any $R$-module $X$ and $g\in G$. The result is used freely in the section.

\begin{prop} Assume that $R$ is a locally bounded $K$-category, $G$ acts freely on the objects of $R$ and $F:R\ra A\cong R\slash G$ is the Galois covering. The functor $F_{\lambda}:\mod(R)\ra\mod(A)$ induces natural isomorphisms $$\nu_{X,Y}:\bigoplus_{g\in G}{}_{R}({}^{g}X,Y)\ra{}_{A}(F_{\lambda}(X),F_{\lambda}(Y))$$ of vector spaces, for any $X,Y\in\mod(R)$, given by $\nu_{X,Y}((f_{g})_{g\in G})=\sum_{g\in G}F_{\lambda}(f_{g})$ where $f_{g}:{}^{g}X\ra Y$, for any $g\in G$.
\end{prop}

We denote by $\Add(\mod(R))$ the full subcategory of $\Mod(R)$ whose objects are arbitrary direct sums of finite dimensional $R$-modules. Assume that an $R$-homomorphism $f:\bigoplus_{j\in J}M_{j}\ra\bigoplus_{i\in I}N_{i}$ in $\Add(\mod(R))$ is defined by homomorphisms $f_{ij}:M_{j}\ra N_{i}$, for $i\in I$, $j\in J$. Observe that for any $j\in J$ we have $f_{ij}\neq 0$ only for finite number of $i\in I$. Indeed, this follows from the fact that $M_{j}$ is finite dimensional, for any $j\in J$.

Assume that $T\in\CF(R)$. We define the functor $\wh{T}:\Add(\mod(R))\ra\Mod(K)$ in the following way. Assume that $\bigoplus_{j\in J}M_{j}$ is an object of $\Add(\mod(R))$ and $f:\bigoplus_{j\in J}M_{j}\ra\bigoplus_{i\in I}N_{i}$ is an $R$-homomorphism in $\Add(\mod(R))$ defined by $f_{ij}:M_{j}\ra N_{i}$, for $i\in I$, $j\in J$. Then we set $\wh{T}(\bigoplus_{j\in J}M_{j})=\bigoplus_{j\in J}T(M_{j})$ and $\wh{T}(f):\bigoplus_{i\in I}T(N_{i})\ra\bigoplus_{j\in J}T(M_{j})$ is defined by $T(f_{ij}):T(N_{i})\ra T(M_{j})$, for $i\in I$, $j\in J$. Observe that $\wh{T}$ equals $T$ on $\mod(R)$.

The functor $\Phi:\CF(R)\ra\CF(A)$ is defined as follows. We set $\Phi(T)=\wh{T}\circ F_{\bullet}$, for any $T\in\CF(R)$. First note that $\Phi(T)\in\CG(A)$, that is, $\Phi(T):\mod(A)\ra\mod(K)$. Indeed, if $X\in\mod(A)$, then $X\cong F_{\lambda}(M)$ for some $M\in\mod(R)$. Thus we get $$\Phi(T)(X)=\Phi(T)(F_{\lambda}(M))=\wh{T}(F_{\bullet}(F_{\lambda}(M)))\cong\wh{T}(\bigoplus_{g\in G}{}^{g}M)=\bigoplus_{g\in G} T({}^{g}M)$$ which is finite dimensional from Lemma 5.1. Assume $\alpha\in{}_{A}(F_{\lambda}(X),F_{\lambda}(Y))$, for some $X,Y\in\mod(R)$, $\alpha=\sum_{g\in G}F_{\lambda}(f_{g})$ where $f_{g}:{}^{g}X\ra Y$, for any $g\in G$. It is worth to note that the homomorphism $$\Phi(T)(\alpha):\bigoplus_{g\in G}T({}^{g}Y)\ra\bigoplus_{g\in G}T({}^{g}X)$$ of vector spaces is defined by homomorphisms $T({}^{g}f_{g^{-1}h}):T({}^{g}Y)\ra T({}^{h}X)$, for any $g,h\in G$. Indeed, this follows from the fact that $F_{\bullet}(\alpha):\bigoplus_{g\in G}{}^{g}X\ra\bigoplus_{g\in G}{}^{g}Y$ is defined by homomorphisms ${}^{g}f_{g^{-1}h}:{}^{h}X\ra {}^{g}Y$, for any $g,h\in G$, see \cite{BoGa}.

Assume that $T_{1},T_{2}\in\CF(R)$ and $\iota:T_{1}\ra T_{2}$ is a homomorphism of functors. Then $\Phi(\iota):\Phi(T_{1})\ra\Phi(T_{2})$ is defined as follows. Assume $X\in\mod(R)$. We set $$\Phi(\iota)_{F_{\lambda}(X)}:\bigoplus_{g\in G}T_{1}({}^{g}X)\ra\bigoplus_{g\in G}T_{2}({}^{g}X)$$ to be the homomorphism of vector spaces defined by $\iota_{{}^{g}X}:T_{1}({}^{g}X)\ra T_{2}({}^{g}X)$, for $g\in G$. It is easy to see that $\Phi(\iota)_{F_{\lambda}(X)}=\Phi(\iota)_{F_{\lambda}({}^{g}X)}$, for any $g\in G$ and $X\in\mod(A)$. Moreover, assume that $\alpha\in{}_{A}(F_{\lambda}(X),F_{\lambda}(Y))$ and $\alpha=\sum_{g\in G}F_{\lambda}(f_{g})$, $f_{g}:{}^{g}X\ra Y$ for $g\in G$. Since $\iota_{{}^{g}X}T_{1}({}^{g}f_{1_{G}})=T_{2}({}^{g}f_{1_{G}})\iota_{{}^{g}Y}$, for any $g\in G$, we easily get the equality $$\Phi(\iota)_{F_{\lambda}(X)}\Phi(T_{1})(\alpha)=\Phi(T_{2})(\alpha)\Phi(\iota)_{F_{\lambda}(Y)}.$$ Hence $\Phi(\iota):\Phi(T_{1})\ra\Phi(T_{2})$ is a homomorphism of functors.

The above definitions give rise to a covariant exact functor $\Phi:\CF(R)\ra\CG(A)$. The following proposition shows that $\Phi(T)=\wh{T}\circ F_{\bullet}\in\CF(A)$, for any $T\in\CF(R)$, and thus $\Phi:\CF(R)\ra\CF(A)$.

\begin{prop} Assume that $M,N\in\mod(R)$ and $f:M\ra N$ is an $R$-homomorphism. 
\begin{enumerate}[\rm(1)]
	\item The homomorphism $\Phi({}_{R}(-,f))$ is isomorphic with the homomorphism $${}_{R}(F_{\bullet}(-),f):{}_{R}(F_{\bullet}(-),M)\ra{}_{R}(F_{\bullet}(-),N).$$
	\item Assume that $T=\Coker{}_{R}(-,f)$. Then $\Phi(T)=\Coker{}_{A}(-,F_{\lambda}(f))\in\CF(A)$.
\end{enumerate}
\end{prop}

{\bf Proof.} (1) Recall that ${}_{R}(F_{\bullet}(F_{\lambda}(X)),Y)\cong\bigoplus_{g\in G}{}_{R}({}^{g}X,Y)$, for any modules $X,Y\in\mod(R)$, and this isomorphism is natural. It is easy to see that $${}_{R}(F_{\bullet}(F_{\lambda}(X)),f):\bigoplus_{g\in G}{}_{R}({}^{g}X,M)\ra\bigoplus_{g\in G}{}_{R}({}^{g}X,N)$$ is defined by homomorphisms ${}_{R}({}^{g}X,f):{}_{R}({}^{g}X,M)\ra{}_{R}({}^{g}X,N)$, for $g\in G$. Since $\Phi({}_{R}(-,f))$ is defined by the same homomorphisms, the assertion follows.

(2) There is an exact sequence of functors ${}_{R}(-,M)\stackrel{{}_{R}(-,f)}{\longrightarrow}{}_{R}(-,N)\ra T\ra 0$ and since the functor $\Phi:\CF(R)\ra\CG(A)$ is exact, the sequence of functors $$\Phi({}_{R}(-,M))\stackrel{\Phi({}_{R}(-,f))}{\longrightarrow}\Phi({}_{R}(-,N))\ra \Phi(T)\ra 0$$ is exact. Recall that $(F_{\bullet},F_{\rho})$ is an adjoint pair and thus the following diagram $$\xymatrix{{}_{R}(F_{\bullet}(-),M)\ar[rr]^{{}_{R}(F_{\bullet}(-),f)}\ar[d]^{\cong}&&{}_{R}(F_{\bullet}(-),N)\ar[d]^{\cong}\\{}_{A}(-,F_{\rho}(M))\ar[rr]^{{}_{A}(-,F_{\rho}(f))}&&{}_{A}(-,F_{\rho}(N))}$$ is commutative. Since $M,N$ are finite dimensional, we get $F_{\rho}(C)=F_{\lambda}(C)$, for $C=M,N,f$ (see \cite{BoGa} or Section 2). Then $(1)$ yields $$\Phi(T)\cong\Coker_{\snull}\Phi({}_{R}(-,f))\cong\Coker{}_{R}(F_{\bullet}(-),f)\cong\Coker{}_{A}(-,F_{\lambda}(f))$$ which shows the assertion. \epv

Summing up, we get a covariant exact functor $\Phi:\CF(R)\ra\CF(A)$ such that $$\Phi(T)=\wh{T}\circ F_{\bullet}\cong\Coker{}_{A}(-,F_{\lambda}(f)),$$ for any $T=\Coker{}_{R}(-,f)\in\CF(R)$. Note that Proposition 5.3 (2) shows in particular that $\Phi({}_{R}(-,M))={}_{A}(-,F_{\lambda}(M))$ and $\Phi({}_{R}(-,f))={}_{A}(-,F_{\lambda}(f))$, for any $R$-modules $M,N$ and $R$-homomorphism $f:M\ra N$. 

Recall that in our situation the functor $F_{\lambda}:\mod(R)\ra\mod(A)$ is a covering of module categories, see Theorem 2.1. The description of $\Phi:\CF(R)\ra\CF(A)$ suggests that it plays the role of the push-down of $F_{\lambda}$, see Section 2 for the precise definition of the push-down functor of a Galois covering of $K$-categories. A natural analogue of the pull-up functor of $F_{\lambda}$ is the functor $\Psi:\CF(A)\ra\CG(R)$ such that $\Psi=(-)\circ F_{\lambda}$. We show in the sequel some relations between $\Phi$ and $\Psi$.

Our aim is to prove some basic properties of the functor $\Phi:\CF(R)\ra\CF(A)$. First observe that the group $G$ acts on $\CF(R)$. Indeed, given a functor $T\in\CG(R)$ and $g\in G$ we define $gT\in\CG(R)$ such that $(gT)(X)=T({}^{g^{-1}}X)$ and $(gT)(f)=T({}^{g^{-1}}f)$, for any module $X\in\mod(R)$ and homomorphism $f\in\mod(R)$. Note that there are isomorphisms $g({}_{R}(-,M))\cong{}_{R}(-,{}^{g}M)$, for any $g\in G$ and $M\in\mod(R)$. Hence $g(\Coker{}_{R}(-,f))\cong\Coker{}_{R}(-,{}^{g}f)$, for any $f\in\mod(R)$ and $g\in G$, and thus the action of $G$ on $\CG(R)$ restricts to $\CF(R)$. We show below that this action is free.

\begin{lm} Assume that $T\in\CF(R)$ and $T\neq 0$. If $g\in G$ and $gT\cong T$, then $g=1$. Thus the group $G$ acts freely on $\CF(R)$.
\end{lm}

{\bf Proof.} Assume that $g\in G$, $gT\cong T$ and $g\neq 1$. Let $X$ be an $R$-module such that $T(X)\neq 0$. Since $gT\cong T$, we get $T(X)\cong T({}^{g^{-n}}X)\neq 0$, for any $n\in\NN$. The elements $g^{-n}$, for $n\geq 0$, are pairwise different, because $G$ is torsion-free. Hence $T({}^{h}X)\neq 0$ for infinitely many $h\in G$. This contradicts Lemma 5.1. Consequently, $gT\cong T$ yields $g=1$, for any $g\in G$, and thus $G$ acts freely on $\CF(R)$. \epv

We define $\Im(\Phi)=\Phi(\CF(R))$ as the full subcategory of $\CF(A)$ whose objects are isomorphic to objects of the form $\Phi(T)$, for $T\in\CF(R)$. 

The following theorem shows that the functor $\Phi:\CF(R)\ra\CF(A)$ induces a Galois covering of additive categories $\CF(R)\ra\Im(\Phi)$ with respect to the action of $G$ introduced above. The properties of $\Phi$ listed in the theorem are applied in Section 6 in proofs of our main results. Additionally, the assertions (2) and (3) exhibit relations between $\Phi$ and $\Psi$. In particular, (3) shows that $\Psi$ is in fact the left adjoint of $\Phi$.

\begin{thm} Assume that $R$ is a locally support-finite $K$-category, $G$ a torsion-free admissible group of $K$-linear automorphisms of $R$ and $F:R\ra A\cong R\slash G$ the Galois covering. Assume that $T,T_{1},T_{2}\in\CF(R)$ and $U\in\CF(A)$. The following assertions hold. 
\begin{enumerate}[\rm(1)]
	\item There exists an isomorphism $\Phi(T)\cong\Phi(gT)$, for any $g\in G$. 
	\item There exists an isomorphism $\Psi(\Phi(T))\cong\bigoplus_{g\in G}gT$. In consequence, if functors $T_{1},T_{2}$ have local endomorphism rings, then $\Phi(T_{1})\cong\Phi(T_{2})$ implies $T_{1}\cong gT_{2}$, for some $g\in G$. 
	\item There exists an isomorphism of vector spaces ${}_{\CG(R)}(\Psi(U),T)\cong {}_{\CG(A)}(U,\Phi(T))$. In particular, the functor $\Phi:\CF(R)\ra\CF(A)$ induces isomorphisms $$\mu_{T_{1},T_{2}}:\bigoplus_{g\in G}{}_{\CF(R)}(gT_{1},T_{2})\ra{}_{\CF(A)}(\Phi(T_{1}),\Phi(T_{2}))$$ of vector spaces given by $\mu_{T_{1},T_{2}}((\iota_{g})_{g\in G})=\sum_{g\in G}\Phi(\iota_{g})$ where $\iota_{g}:gT_{1}\ra T_{2}$ is a homomorphism of functors.
	\item There exists is an epimorphism $\epsilon:\Phi(T')\ra U$, for some $T'\in\CF(R)$.
\end{enumerate}
\end{thm}

{\bf Proof.} (1) Assume that $g\in G$ and $T=\Coker{}_{R}(-,f)$, for some $R$-homomorphism $f:M\ra N$. Then $gT\cong\Coker{}_{R}(-,{}^{g}f)$, $F_{\lambda}(f)=F_{\lambda}({}^{g}f)$ and thus $$\Phi(gT)\cong\Coker{}_{A}(-,F_{\lambda}({}^{g}f))\cong\Coker{}_{A}(-,F_{\lambda}(f))=\Phi(T).$$

(2) Observe that there are isomorphisms $$\Psi(\Phi(T))(X)=\wh{T}(F_{\bullet}(F_{\lambda}(X)))=\wh{T}(\bigoplus_{g\in G}{}^{g}X)\cong\bigoplus_{g\in G} T({}^{g}X)\cong(\bigoplus_{g\in G}gT)(X),$$ for any $X\in\mod(R)$. These isomorphisms are natural which follows directly from the definition of $\wh{T}:\Add(\mod(R))\ra\Mod(K)$. This yields $\Psi(\Phi(T))\cong\bigoplus_{g\in G}gT$.

Assume that $\Phi(T_{1})\cong\Phi(T_{2})$. Then $\Psi(\Phi(T_{1}))\cong\Psi(\Phi(T_{2}))$ and hence there is an isomorphism $\bigoplus_{g\in G}gT_{1}\cong\bigoplus_{g\in G}gT_{2}$. Note that $gT_{1},gT_{2}$ have local endomorphism rings, for any $g\in G$. Hence the Krull-Remak-Schmidt theorem (see for example Appendix E of \cite{Pr2}) yields $T_{1}\cong gT_{2}$, for some $g\in G$.

(3) We show that ${}_{\CG(R)}(U\circ F_{\lambda},T)\cong {}_{\CG(A)}(U,\wh{T}\circ F_{\bullet})$. Assume that $\iota:U\circ F_{\lambda}\ra T$ is a homomorphism of functors. This means that $T(f)\iota_{Y}=\iota_{X}U(F_{\lambda}(f))$, for any $R$-module homomorphism $f:X\ra Y$. We define $\tau:U\ra\wh{T}\circ F_{\bullet}$ in the following way. Assume that $p_{{}^{g}X}:\bigoplus_{g\in G}T({}^{g}X)\ra T({}^{g}X)$ is the split epimorphism, for any $R$-module $X$ and $g\in G$. Recall that $(\wh{T}\circ F_{\bullet})(F_{\lambda}(X))\cong\bigoplus_{g\in G}T({}^{g}X)$ and let $$\tau_{F_{\lambda}(X)}:U(F_{\lambda}(X))\ra\bigoplus_{g\in G}T({}^{g}X)$$ be defined by homomorphisms $\iota_{{}^{g}X}:U(F_{\lambda}({}^{g}X))\ra T({}^{g}X)$, for any $g\in G$, which means that $p_{{}^{g}X}\tau_{F_{\lambda}(X)}=\iota_{{}^{g}X}$. The definition is correct since $F_{\lambda}({}^{g}X)\cong F_{\lambda}(X)$ and $\tau_{F_{\lambda}(X)}=\tau_{F_{\lambda}({}^{g}X)}$, for any $g\in G$. We show that the above vector space homomorphisms define a homomorphism of functors $\tau:U\ra\wh{T}\circ F_{\bullet}$. Indeed, assume that $\alpha\in{}_{A}(F_{\lambda}(X),F_{\lambda}(Y))$ and $\alpha=\sum_{g\in G}F_{\lambda}(f_{g})$ where $f_{g}:{}^{g}X\ra Y$, for any $g\in G$. Recall that $$\wh{T}(F_{\bullet}(\alpha)):\bigoplus_{g\in G}T({}^{g}Y)\ra\bigoplus_{g\in G}T({}^{g}X)$$ is defined by homomorphisms $T({}^{g}f_{g^{-1}h}):T({}^{g}Y)\ra T({}^{h}X)$, for any $g,h\in G$. This yields $$p_{{}^{h}X}\wh{T}(F_{\bullet}(\alpha))=\sum_{g\in G}T({}^{hg^{-1}}f_{g})p_{{}^{hg^{-1}}Y},$$ for any $h\in G$. Therefore we get the following equalities $$p_{{}^{h}X}\tau_{F_{\lambda}(X)}U(\alpha)=\iota_{{}^{h}X}\sum_{g\in G}U(F_{\lambda}(f_{g}))=\sum_{g\in G}\iota_{{}^{h}X}U(F_{\lambda}({}^{hg^{-1}}f_{g}))=$$$$=\sum_{g\in G}T({}^{hg^{-1}}f_{g})\iota_{{}^{hg^{-1}}Y}=\sum_{g\in G}T({}^{hg^{-1}}f_{g})p_{{}^{hg^{-1}}Y}\tau_{F_{\lambda}(Y)}=p_{{}^{h}X}\wh{T}(F_{\bullet}(\alpha))\tau_{F_{\lambda}(Y)},$$ for any $h\in G$. This yields $\tau_{F_{\lambda}(X)}U(\alpha)=\wh{T}(F_{\bullet}(\alpha))\tau_{F_{\lambda}(Y)}$ and hence $\tau:U\ra\wh{T}\circ F_{\bullet}$ is a homomorphism of functors. It is easy to see that the map $\iota\mapsto\tau$ defines a homomorphism of the appropriate vector spaces. We denote this homomorphism by $\chi$, that is, $\chi(\iota)=\tau$.

Similar arguments imply that if $\tau:U\ra\wh{T}\circ F_{\bullet}$ is a homomorphism of functors, then the vector space homomorphisms $\iota_{{}^{g}X}=p_{{}^{g}X}\tau_{F_{\lambda}(X)}$, for $X\in\mod(R)$ and $g\in G$, define a homomorphism of functors $\iota:U\circ F_{\lambda}\ra T$. It is easy to see that the map $\tau\mapsto\iota$ defines a homomorphism the appropriate vector spaces. Since the above vector space homomorphisms are mutually inverse, we get ${}_{\CG(A)}(U,\wh{T}\circ F_{\bullet})\cong{}_{\CG(R)}(U\circ F_{\lambda},T)$.

The second part of the assertion follows from the first one by putting $U=\wh{T_{1}}\circ F_{\bullet}$ and $T=T_{2}$. To be more specific, assume that $\iota_{g}:gT_{1}\ra T_{2}$ is a homomorphism of functors, for any $g\in G$. Note that there is only finite number of $g\in G$ such that $\iota_{g}\neq 0$, because $(\iota_{g})_{g\in G}$ is an element of a direct sum. The homomorphisms $\iota_{g}$, $g\in G$, yield a homomorphism of functors $\iota:(\wh{T_{1}}\circ F_{\bullet})\circ F_{\lambda}\ra T_{2}$ such that the homomorphism $\iota_{X}:\bigoplus_{g\in G}T_{1}({}^{g}X)\ra T_{2}(X)$ is defined by $(\iota_{g})_{X}:T_{1}({}^{g^{-1}}X)\ra T_{2}(X)$, for $g\in G$. The map $(\iota_{g})_{g\in G}\mapsto\iota$ is an isomorphism of vector spaces. Moreover, it is easy to see that $\chi(\iota)=\sum_{g\in G}\Phi(\iota_{g})$ and the assertion follows.

(4) Assume that $U\in\CF(A)$. Then $U$ is finitely generated, thus there is an epimorphism of functors of the form $\epsilon:{}_{A}(-,M)\ra U$ , for some $A$-module $M$. Since $F_{\lambda}$ is dense, we get $M\cong F_{\lambda}(X)$, for some $R$-module $X$, so ${}_{A}(-,M)\cong{}_{A}(-,F_{\lambda}(X))=\Phi({}_{R}(-,X))$. This shows the assertion. \epv

Observe that the assertions (1), (2) and (3) of Theorem 5.5 show that the functor $\Phi:\CF(R)\ra\CF(A)$ induces a Galois covering $\CF(R)\ra\Im(\Phi)$ of additive categories. The assertion (4) does not imply the density of $\Phi$.

\begin{rem} The density of the functor $\Phi$ is hard to verify. However, we guess that $\Phi$ may not be dense without some additional assumptions. Indeed, the functor $\Phi$ plays the role of the push-down functor of $F_{\lambda}:\mod(R)\ra\mod(A)$ and push-down functors are not dense in general (note that $F_{\lambda}$ is dense, because $R$ is a locally support-finite $K$-category, see Section 2). 

The main results of the paper yield $\KG(R)=\KG(A)$. Therefore it is natural to ask whether the condition $\KG(R)=\KG(A)<\infty$ is sufficient for the density of $\Phi$. This is an interesting open problem which is left for further research.

Assume that $U\in\CF(A)$. We show below another example of epimorphism of functors $\epsilon:\Phi(T)\ra U$, for some $T\in\CF(R)$. It is hard to verify whether $\epsilon$ is an isomorphism, but this example is natural to consider.

Recall that $F_{\lambda}:\mod(R)\ra\mod(A)$ is dense, so $U=\Coker{}_{A}(-,\alpha)$, for some $A$-homomorphism $\alpha:F_{\lambda}(M)\ra F_{\lambda}(N)$. Proposition 5.2 yields $F_{\lambda}$ induces an isomorphism $${}_{A}(F_{\lambda}(M),F_{\lambda}(N))\cong\bigoplus_{g\in G}{}_{R}({}^{g}M,N).$$ Hence assume that $\alpha=\sum_{g\in G}F_{\lambda}(f_{g})$, $f_{g}:{}^{g}M\ra N$ for $g\in G$, and $g_{1},...,g_{n}$ are the only elements of $G$ such that $f_{g_{i}}\neq 0$, for $i=1,...,n$. Set $N=\bigoplus_{i=1}^{n}{}^{g_{i}^{-1}}N$, let $\ov{f}:M\ra N$ be the $R$-homomorphism such that $$\ov{f}=\left[\begin{matrix}{}^{g_{1}^{-1}}f_{g_{1}}&{}^{g_{2}^{-1}}f_{g_{2}}&\hdots& {}^{g_{n}^{-1}}f_{g_{n}}\end{matrix}\right]^{t}$$ and define $T=\Coker{}_{R}(-,\ov{f})$. We show that there is an epimorphism of functors $\epsilon:\Phi(T)\ra U$. 

Recall that $F_{\lambda}({}^{g_{i}^{-1}}N)\cong F_{\lambda}(N)$ and $F_{\lambda}({}^{g_{i}^{-1}}f_{g_{i}})=F(f_{g_{i}})$, for any $i=1,...,n$. Let $\pi:\bigoplus_{i=1}^{n}F_{\lambda}(N)\ra F_{\lambda}(N)$ be the epimorphism given by the matrix $\left[\begin{matrix}1&1&\hdots&1\end{matrix}\right]$. Then we get $\alpha=\sum_{i=1}^{n}F_{\lambda}(f_{g_{i}})=\pi F_{\lambda}(\ov{f})$, so ${}_{A}(-,\alpha)={}_{A}(-,\pi){}_{A}(-,F_{\lambda}(\ov{f}))$. It is easy to see that ${}_{A}(-,\pi):{}_{A}(-,\bigoplus_{i=1}^{n}F_{\lambda}(N))\ra{}_{A}(-, F_{\lambda}(N))$ is an epimorphism and hence ${}_{A}(-,\pi)$ induces an epimorphism $\epsilon:\Coker{}_{A}(-,F_{\lambda}(\ov{f}))\ra\Coker{}_{A}(-,\alpha)$. This shows the claim since $\Phi(T)=\Coker{}_{A}(-,F_{\lambda}(\ov{f}))$. 

Observe that $\epsilon:\Phi(T)\ra U$ is an isomorphism if and only if for any $R$-module $X$ and any $A$-homomorphisms $\gamma_{1},...,\gamma_{n}:F_{\lambda}(X)\ra F_{\lambda}(N)$ the following condition holds: $\gamma_{1}+...+\gamma_{n}=\alpha(F_{\lambda}(f_{g_{1}})+...+F_{\lambda}(f_{g_{n}}))$, for some $\alpha:F_{\lambda}(X)\ra F_{\lambda}(M)$, if and only if $\gamma_{i}=\beta F_{\lambda}(f_{g_{i}})$, for some $\beta:F_{\lambda}(X)\ra F_{\lambda}(M)$, for any $i=1,...,n$. It is hard to verify whether this holds in general. \epv
\end{rem}   

\section{Krull-Gabriel dimension of $R$ and $A$}

Throughout we assume that $R$ is a locally support-finite $K$-category, $G$ a torsion-free admissible group of $K$-linear automorphisms of $R$ and $F:R\ra A\cong R\slash G$ the Galois covering. In particular, there exists the functor $\Phi:\CF(R)\ra\CF(A)$ which induces a Galois covering $\CF(R)\ra\Im(\Phi)$ of additive categories, see Theorem 5.5. Assume that $(\CF(R)_{\alpha})_{\alpha}$ and $(\CF(A)_{\alpha})_{\alpha}$ are the Krull-Gabriel filtrations of the categories $\CF(R)$ and $\CF(A)$, respectively. This section is devoted to showing that $T\in\CF(R)_{\alpha}$ if and only if $\Phi(T)\in\CF(A)_{\alpha}$, for any $T\in\CF(R)$ and ordinal number $\alpha$. Then it follows from Lemma 3.2 (3) and Theorem 5.5 (4) that $\KG(R)=\KG(A)$. Moreover, we give a sufficient condition for the existence of isomorphisms $\CF(R)_{\alpha}\slash G\cong\CF(A)_{\alpha}$, for any ordinal $\alpha$.

Assume that $B$ is $R$ or $A$. We recall from Section 3 that a functor $S\in\CF(B)$ is simple in the category $\CF(B)\slash\CF(B)_{\alpha}$ if and only if $S$ is $\CF(B)_{\alpha}$-simple, see Lemma 3.1 for the details. We apply this fact freely throughout the section.


We show that the Krull-Gabriel filtration of the category $\CF(R)$ is $G$-invariant. This means that $T\in\CF(R)_{\alpha}$ yields $gT\in\CF(R)_{\alpha}$, for any ordinal $\alpha$ and $g\in G$.

\begin{prop} Assume that $R$ is a locally support-finite $K$-category and $G$ a torsion-free admissible group of $K$-linear automorphisms of $R$.
\begin{enumerate}[\rm(1)]
	\item Assume that $\alpha$ is an ordinal and the functor $S\in\CF(R)$ is simple in the category $\CF(R)_{\alpha}$. Then the functor $gS\in\CF(R)$ is simple in the category $\CF(R)_{\alpha}$, for any $g\in G$.
	\item If $T\in\CF(R)_{\alpha}$, then $gT\in\CF(R)_{\alpha}$, for any ordinal $\alpha$ and $g\in G$. Thus the Krull-Gabriel filtration of the category $\CF(R)$ is $G$-invariant.
	\item The group $G$ acts freely on $\CF(R)_{\alpha}$, for any ordinal $\alpha$.
\end{enumerate}
\end{prop}

{\bf Proof.} We show (1) and (2) simultaneously. We proceed by transfinite induction with respect to $\alpha$. Assume that $\alpha=0$ and $S$ is a simple functor in $\CF(R)_{0}$, which means that $S$ is simple in $\CF(R)$. Assume that $g\in G$. We show that $gS$ is also simple in $\CF(R)$. Indeed, we have $S\cong{}_{R}(-,M)\slash\rad_{R}(-,M)$ for some indecomposable $R$-module $M$, see Section 4. It is easy to see that $g(\rad_{R}(-,M))\cong\rad_{R}(-,{}^{g}M)$ and hence $$gS\cong g(\frac{{}_{R}(-,M)}{\rad_{R}(-,M)})\cong \frac{g({}_{R}(-,M))}{g(\rad_{R}(-,M))}\cong\frac{{}_{R}(-,{}^{g}M)}{\rad_{R}(-,{}^{g}M)}.$$ Since ${}^{g}M$ is indecomposable (because $G$ acts freely on $\ind(R)$), this gives that $gS$ is simple in $\CF(R)$ and thus it is simple in $\CF(R)_{0}$. 

Assume that $T\in\CF(R)_{0}$, that is, $T$ has finite length in $\CF(R)$. This means that there exists a chain $0=T_{0}\subseteq T_{1}\subseteq ...\subseteq T_{t}=T$ of subobjects of $T$ such that $T_{i+1}\slash T_{i}$ is simple in $\CF(R)$, for any $i=0,...,t-1$. Therefore we get a chain $0=gT_{0}\subseteq gT_{1}\subseteq ...\subseteq gT_{t}=gT$ of subobjects of $gT$, for any $g\in G$, such that $gT_{i+1}\slash gT_{i}\cong  g(T_{i+1}\slash T_{i})$ is simple in $\CF(R)$. Hence $gT$ has finite length in $\CF(R)$. This implies that $\CF(R)_{0}$ is $G$-invariant.


Assume that $\alpha$ is an ordinal and the category $\CF(R)_{\alpha}$ is $G$-invariant. We show that the category $\CF(R)_{\alpha+1}$ is $G$-invariant. Assume that $S$ is simple in $\CF(R)\slash\CF(R)_{\alpha}$. We show that $gS$ is also simple in $\CF(R)\slash\CF(R)_{\alpha}$, for any $g\in G$. Indeed, the functor $gS$ is not an element of $\CF(R)_{\alpha}$ since otherwise $g^{-1}(gS)\cong S$ belongs to $\CF(R)_{\alpha}$, a contradiction. Assume that $P$ is a subobject of $gS$. It follows that $g^{-1}P$ is a subobject of $S$ and hence $g^{-1}P\in\CF(R)_{\alpha}$ or $(S\slash g^{-1}P)\in\CF(R)_{\alpha}$. Then $g(g^{-1}P)\cong P\in\CF(R)_{\alpha}$ or $g(S\slash g^{-1}P)\cong (gS\slash P)\in\CF(R)_{\alpha}$, because $\CF(R)_{\alpha}$ is $G$-invariant. This implies that $gS$ is simple in $\CF(R)\slash\CF(R)_{\alpha}$. 

Assume that $T$ has finite length in $\CF(R)\slash\CF(R)_{\alpha}$ (so $T\in\CF(R)_{\alpha+1}$). Then there is a chain $0=T_{0}\subseteq T_{1}\subseteq ...\subseteq T_{t}=T$ of subobjects of $T$ such that $S_{i}=T_{i+1}\slash T_{i}$ is simple in $\CF(R)\slash\CF(R)_{\alpha}$, for any $i=0,...,t-1$. This implies that $gS_{i}=gT_{i+1}\slash gT_{i}$ is simple in $\CF(R)\slash\CF(R)_{\alpha}$ and thus the chain $0=gT_{0}\subseteq gT_{1}\subseteq ...\subseteq gT_{t}=gT$ of subobjects of $gT$ is a composition series of $gT$. Consequently, $gT\in\CF(R)_{\alpha+1}$ and thus the category $\CF(R)_{\alpha+1}$ is $G$-invariant. 

If $\beta$ is a limit ordinal, then $\CF(R)_{\beta}=\bigcup_{\alpha<\beta}\CF(R)_{\alpha}$. Hence $\CF(R)_{\beta}$ is $G$-invariant, if we assume that all $\CF(R)_{\alpha}$ are. 

The above arguments show (1) and (2) (by transfinite induction). The assertion of (3) follows from (1), (2) and Lemma 5.4. \epv

A monomorphism of functors is \textit{proper} if and only if it is not an isomorphism. The following fact is applied in the sequel.

\begin{lm} Assume that $T_{1},T_{2}\in\CF(R)$ and $\iota:T_{1}\ra T_{2}$ is a proper monomorphism of functors. Then $\Phi(\iota):\Phi(T_{1})\ra\Phi(T_{2})$ is a proper monomorphism.
\end{lm}

{\bf Proof.} Assume that $\iota:T_{1}\ra T_{2}$ is a proper monomorphism. Then the homomorphism $\Phi(\iota):\Phi(T_{1})\ra\Phi(T_{2})$ is a monomorphism, because $\Phi$ is exact. We show that $\Phi$ is proper. Indeed, if $\iota$ is proper, then there is an $R$-module $X$ such that $\dim_{K}(T_{1}(X))<\dim_{K}(T_{2}(X))$. Since $\dim_{K}(T_{1}({}^{g}X))\leq\dim_{K}(T_{2}({}^{g}X))$, for any $g\in G$, and $T_{1}({}^{g}X),T_{2}({}^{g}X)\neq 0$ only for finite number of $g\in G$ (see Lemma 5.1), we get $$\dim_{K}(\bigoplus_{g\in G}T_{1}({}^{g}X))<\dim_{K}(\bigoplus_{g\in G}T_{2}({}^{g}X)).$$ This yields $\Phi(\iota):\Phi(T_{1})\ra\Phi(T_{2})$ is proper, because $\Phi(T_{i})(F_{\lambda}(X))\cong\bigoplus_{g\in G}T_{i}({}^{g}X)$, for $i=1,2$, see Section 5. \epv

The following theorem and corollary are the main results of the paper. In their proofs we apply the fact that the functor $\Phi:\CF(R)\ra\CF(A)$ induces a Galois covering $\CF(R)\ra\Im(\Phi)$ of additive categories, see Theorem 5.5.

\begin{thm} Assume that $R$ is a locally support-finite $K$-category, $G$ a torsion-free admissible group of $K$-linear automorphisms of $R$ and $F:R\ra A\cong R\slash G$ the Galois covering. The following assertions hold. 
\begin{enumerate}[\rm(1)]
	\item Assume that $T\in\CF(R)$. Then $T\in\CF(R)_{\alpha}$ if and only if $\Phi(T)\in\CF(A)_{\alpha}$, for any ordinal number $\alpha$.
	\item We have $\KG(R)=\KG(A)$.
\end{enumerate} 
\end{thm}

{\bf Proof.} (1) We proceed by transfinite induction with respect to $\alpha$. Assume that $\alpha=0$. The categories $\CF(R)_{0}$ and $\CF(A)_{0}$ consist of all functors of finite length in $\CF(R)$ and $\CF(A)$, respectively. First we show that $S\in\CF(R)$ is simple if and only if $\Phi(S)\in\CF(A)$ is simple. 

Assume that $S\in\CF(R)$ is a simple functor. Then there exists an indecomposable $R$-module $N$ and right minimal almost split $R$-homomorphism $f:M\ra N$ such that $S\cong S^{N}$ and the sequence of functors ${}_{R}(-,M)\stackrel{{}_{R}(-,f)}{\longrightarrow}{}_{R}(-,N)\stackrel{\pi^{N}}{\ra}S^{N}\ra 0$ is a minimal projective presentation of $S^{N}$ in $\CF(R)$, see Section 4. Since $F_{\lambda}$ preserves right minimal almost split homomorphisms (see \cite{DoSk} or Section 2), we get that $${}_{A}(-,F_{\lambda}(M))\stackrel{{}_{A}(-,F_{\lambda}(f))}{\longrightarrow}{}_{A}(-,F_{\lambda}(N))\stackrel{F_{\lambda}(\pi^{N})}{\longrightarrow}\Phi(S^{N})\ra 0$$ is a minimal projective presentation of $\Phi(S^{N})$. This yields that $\Phi(S^{N})\cong\Phi(S)$ is a simple functor in $\CF(A)$.


Assume that $\Phi(S)\in\CF(A)$ is a simple functor. If $S\in\CF(R)$ is not simple, there is $S'\in\CF(R)$, $S'\neq 0$ and a proper monomorphism $\iota:S'\ra S$. Then $\Phi(S')\neq 0$ and $\Phi(\iota):\Phi(S')\ra\Phi(S)$ is a proper monomorphism, see Lemma 6.2. This contradicts the fact that $\Phi(S)$ is simple. Thus if $\Phi(S)\in\CF(A)$ is simple, then $S$ is simple. 



Consequently, $S$ is simple in $\CF(R)$ if and only if $\Phi(S)$ is simple in $\CF(A)$. This implies that the functor $\Phi$ preserves composition series, because $\Phi$ is exact and $T_{1}\subseteq T_{2}$ yields $\Phi(T_{1})\subseteq\Phi(T_{2})$. Therefore, if $T$ has finite length, then $\Phi(T)$ has finite length, for any $T\in\CF(R)$. 

For the converse implication, observe that if $T$ does not have finite length, there are arbitrarily long chains $0=T_{0}\subseteq T_{1}\subseteq ...\subseteq T_{n}=T$ of subobjects of $T$ such that $T_{i}\subseteq T_{i+1}$ is a proper monomorphism, for any $i=0,...,n-1$. Hence Lemma 6.2 implies that there are arbitrarily long chains of subobjects of $\Phi(T)$. This yields $\Phi(T)$ does not have finite length.

It follows that $T$ has finite length if and only if $\Phi(T)$ has finite length, for any $T\in\CF(R)$. Hence $T\in\CF(R)_{0}$ if and only if $\Phi(T)\in\CF(A)_{0}$, for any $T\in\CF(R)$.

Assume that $\alpha$ is an ordinal and $T\in\CF(R)_{\alpha}$ if and only if $\Phi(T)\in\CF(A)_{\alpha}$, for any $T\in\CF(R)$. We show that $T\in\CF(R)_{\alpha+1}$ if and only if $\Phi(T)\in\CF(A)_{\alpha+1}$, for any $T\in\CF(R)$. First we show that $S$ is simple in the quotient category $\CF(R)\slash\CF(R)_{\alpha}$ if and only if $\Phi(S)$ is simple in the quotient category $\CF(A)\slash\CF(A)_{\alpha}$. 

Assume that $S$ is simple in $\CF(R)\slash\CF(R)_{\alpha}$, that is, $S\notin\CF(R)_{\alpha}$ and $S'\in\CF(R)_{\alpha}$ or $S\slash S'\in\CF(R)_{\alpha}$, for any $S'\subseteq S$. Since $S\notin\CF(R)_{\alpha}$, we get $\Phi(S)\notin\CF(A)_{\alpha}$. Assume that $P\subseteq\Phi(S)$. We show that $P\in\CF(A)_{\alpha}$ or $\Phi(S)\slash P\in\CF(A)_{\alpha}$ which yields $\Phi(S)$ is simple in $\CF(A)\slash\CF(A)_{\alpha}$. Assume that the inclusion $P\subseteq\Phi(S)$ is given by the monomorphism of functors $\iota:P\ra\Phi(S)$. Theorem 5.5 (4) implies that there is an epimorphism of functors $\epsilon:\Phi(T)\ra P$, for some $T\in\CF(R)$. Set $\alpha=\iota\epsilon:\Phi(T)\ra\Phi(S)$ and note that $\Im(\alpha)\cong P$. Theorem 5.5 (3) implies that there are elements $g_{1},...,g_{m}\in G$ and homomorphisms of functors $\alpha_{i}:g_{i}T\ra S$, for $i=1,...,m$, such that $\alpha=\sum_{i=1}^{m}\Phi(\alpha_{i})$. We set $Q=\sum_{i=1}^{m}\Im(\alpha_{i})\subseteq S$. Since $\Phi$ is exact, we get $\Im(\Phi(\alpha_{i}))=\Phi(\Im(\alpha_{i}))$, for any $i=1,...,m$, and hence $$\Phi(Q)=\Phi(\sum_{i=1}^{m}\Im(\alpha_{i}))=\sum_{i=1}^{m}\Phi(\Im(\alpha_{i}))=\sum_{i=1}^{m}\Im(\Phi(\alpha_{i}))=$$$$=\Im(\sum_{i=1}^{m}\Phi(\alpha_{i}))=\Im(\alpha)\cong P.$$ Since $Q\subseteq S$ and $S$ is simple in $\CF(R)\slash\CF(R)_{\alpha}$, we have $Q\in\CF(R)_{\alpha}$ or $S\slash Q\in\CF(R)_{\alpha}$. If $Q\in\CF(R)_{\alpha}$, then $P\cong\Phi(Q)\in\CF(A)_{\alpha}$. If $S\slash Q\in\CF(R)_{\alpha}$, then $$\Phi(S)\slash P\cong\Phi(S)\slash\Phi(Q)=\Phi(S\slash Q)\in\CF(A)_{\alpha}.$$ These arguments imply that $\Phi(S)$ is simple in $\CF(R)\slash\CF(R)_{\alpha}$. 

Assume that $\Phi(S)$ is simple in $\CF(A)\slash\CF(A)_{\alpha}$. Then $S\notin\CF(R)_{\alpha}$, because otherwise $\Phi(S)\in\CF(A)_{\alpha}$. If $T\subseteq S$, then $\Phi(T)\subseteq\Phi(S)$ and hence we have $\Phi(T)\in\CF(A)_{\alpha}$ or $\Phi(S\slash T)=\Phi(S)\slash\Phi(T)\in\CF(A)_{\alpha}$. Then we get $T\in\CF(R)_{\alpha}$ or $S\slash T\in\CF(R)_{\alpha}$ and thus $S$ is simple in $\CF(R)\slash\CF(R)_{\alpha}$. 

Consequently, the functor $S$ is simple in $\CF(R)\slash\CF(R)_{\alpha}$ if and only if the functor $\Phi(S)$ is simple in $\CF(R)\slash\CF(R)_{\alpha}$. This fact implies that the functor $\Phi$ preserves composition series in $\CF(R)\slash\CF(R)_{\alpha}$. Therefore, if $T$ has finite length in $\CF(R)\slash\CF(R)_{\alpha}$, then $\Phi(T)$ has finite length in $\CF(A)\slash\CF(A)_{\alpha}$, for any $T\in\CF(R)$. 

For the converse implication, observe that if $T$ does not have finite length in $\CF(A)\slash\CF(A)_{\alpha}$, there are arbitrarily long chains $0=T_{0}\subseteq T_{1}\subseteq ...\subseteq T_{n}=T$ of subobjects of $T$ such that $T_{i+1}\slash T_{i}$ is not simple in $\CF(R)\slash\CF(R)_{\alpha}$, for some $i=0,...,n-1$. Hence there are arbitrarily long chains $$0=\Phi(T_{0})\subseteq\Phi(T_{1})\subseteq ...\subseteq\Phi(T_{n})=\Phi(T)$$ of subobjects of $\Phi(T)$ such that $\Phi(T_{i+1})\slash\Phi(T_{i})$ is not simple in $\CF(A)\slash\CF(A)_{\alpha}$, for some $i=1,...,n-1$. This yields $\Phi(T)$ does not have finite length in $\CF(A)\slash\CF(A)_{\alpha}$.

It follows that $T$ has finite length in $\CF(R)\slash\CF(R)_{\alpha}$ if and only if $\Phi(T)$ has finite length in $\CF(A)\slash\CF(A)_{\alpha}$, for any $T\in\CF(R)$. Hence $T\in\CF(R)_{\alpha+1}$ if and only if $\Phi(T)\in\CF(A)_{\alpha+1}$, for any $T\in\CF(R)$. Assume that $\beta$ is a limit ordinal and $T\in\CF(R)_{\alpha}$ if and only if $\Phi(T)\in\CF(A)_{\alpha}$, for any ordinal $\alpha<\beta$. Since $\CF(X)_{\beta}=\bigcup_{\alpha<\beta}\CF(X)_{\alpha}$ where $X$ is $R$ or $A$, we get $T\in\CF(R)_{\beta}$ if and only if $\Phi(T)\in\CF(A)_{\beta}$. The assertion of (1) follows by transfinite induction.

(2) The assertion of (2) follows directly from (1), Theorem 5.5 (4) and Lemma 3.2 (3). \epv


\begin{cor} Assume that $R$ is a locally support-finite $K$-category, $G$ a torsion-free admissible group of $K$-linear automorphisms of $R$ and $F:R\ra A$ the Galois covering. If the category $\Im(\Phi)=\Phi(\CF(R))$ is a Serre subcategory of $\CF(A)$, then the functor $\Phi:\CF(R)\ra\CF(A)$ is dense. If this is the case, then $\Phi(\CF(R)_{\alpha})=\CF(A)_{\alpha}$ and $\CF(R)_{\alpha}\slash G\cong\CF(A)_{\alpha}$, for any ordinal number $\alpha$.
\end{cor}

{\bf Proof.} Assume that $\Im(\Phi)$ is a Serre subcategory of $\CF(A)$. If $U\in\CF(A)$, then there is an epimorphism of functors $\epsilon:\Phi(T)\ra U$, for some $T\in\CF(R)$, see Theorem 5.5 (4). Hence $U\in\Im(\Phi)$, because $\Im(\Phi)$ is closed under images of epimorphisms. This implies that the functor $\Phi:\CF(R)\ra\CF(A)$ is dense.

Assume that $\Phi:\CF(R)\ra\CF(A)$ is dense and $\alpha$ is an ordinal number. Assume that $U\in\CF(A)_{\alpha}$. Then $U\cong\Phi(T)$ for some $T\in\CF(R)_{\alpha}$, from Theorem 6.3 (1). Moreover, Theorem 6.3 (1) implies that $T\in\CF(R)_{\alpha}$ yields $\Phi(T)\in\CF(A)_{\alpha}$, for any $T\in\CF(R)_{\alpha}$. Consequently, we get $\Phi(\CF(R)_{\alpha})=\CF(A)_{\alpha}$ and thus Proposition 6.1 (3) and assertions (1), (2), (3) of Theorem 5.5 imply that $\CF(R)_{\alpha}\slash G\cong\CF(A)_{\alpha}$. \epv

\section{Krull-Gabriel dimension of locally support-finite repetitive $K$-categories}

The main result of this section determines the Krull-Gabriel dimension of all locally support-finite repetitive $K$-categories. The proof is based on Lemma 4.3 and two general theorems. The first one is the classification of tame locally support-finite $K$-categories, given by I. Assem and A. Skowro{\'n}ski in \cite{AsSk4}. The second one is the classification of cycle-finite algebras with finite Krull-Gabriel dimension, given by A. Skowro{\'n}ski in \cite{Sk5}.

Recall that it follows from \cite{DoSk2} that a locally bounded $K$-category $R$ is tame if and only if any finite full subcategory of $R$ is tame. Moreover, $R$ is wild if and only if there exists a finite full subcategory of $R$ which is wild.

Assume that $A$ is a $K$-algebra and $1_{A}=e_{1}+...+e_{s}$ is a decomposition of the identity $1_{A}$ of $A$ into sum of orthogonal primitive idempotents. Following \cite{HW}, the \textit{repetitive category} of $A$ is the category $\widehat{A}$ whose objects are $e_{m,i}$, for $m\in\ZZ$, $i=1,...,s$, and the morphism spaces are defined in the following way $$\widehat{A}(e_{m,i},e_{r,j})=\left\{\begin{array}{ccc}e_{j}Ae_{i}, & r=m,\\D(e_{i}Ae_{j}), & r=m+1,\\ 0,& \textnormal{otherwise}\end{array}\right.$$ where $D=\Hom_{K}(-,K)$ denotes the standard duality. Assume that $n\in\NN$. We denote by $\widehat{A}_{n}$ the full subcategory of $\widehat{A}$ whose objects are $e_{m,i}$, for $m=\{-n,...,0,...,n\}$, $i=1,...,s$. Observe that $\widehat{A}_{n}$ is a convex subcategory of $\widehat{A}$. Moreover, $\widehat{A}_{0}\cong A$.

The following theorem is a part of Theorem (B) proved in \cite{AsSk4}.

\begin{thm} Assume that $K$ is an algebraically closed field and $A$ is a finite dimensional basic and connected $K$-algebra. The following conditions are equivalent.
\begin{enumerate}[\rm(1)]
	\item The repetitive $K$-category $\widehat{A}$ is locally support-finite and tame.
	\item There exists an algebra $B$ such that $\widehat{A}\cong\widehat{B}$ and $B$ is tilted of Dynkin type, tilted of Euclidean type or tubular.
\end{enumerate}
\end{thm}

We refer the reader to \cite{AsSiSk} and \cite{Ri} for more information on tilted algebras and tubular algebras, respectively. Assume that $A$ is an algebra. A \textit{cycle} in $\ind(A)$ is a sequence $$M_{0}\stackrel{f_{1}}{\ra} M_{1}\ra\hdots\ra M_{r-1}\stackrel{f_{r}}{\ra}M_{r}=M_{0}$$ of nonzero nonisomorphisms in $\ind(A)$. This cycle is \textit{finite} if and only if the homomorphisms $f_{1},...,f_{r}$ do not belong to $\rad_{A}^{\infty}$. Following \cite{AsSk1}, \cite{AsSk2} we call the algebra $A$ \textit{cycle-finite} if and only if all cycles in $\ind(A)$ are finite. 

The following theorem is a combination of Theorem 1.2 and Corollary 1.5 of \cite{Sk5}.

\begin{thm} Assume that $K$ is an algebraically closed field and $A$ is a finite dimensional basic and connected $K$-algebra. If $A$ is cycle-finite, then $A$ is domestic if and only if $\KG(A)=2$.
\end{thm}

The following theorem is the main result of the section. We apply this fact in Section 8 in order to determine the Krull-Gabriel dimension of standard selfinjective algebras of polynomial growth.

\begin{thm} Assume that $K$ is an algebraically closed field and $A$ is a finite dimensional basic and connected $K$-algebra such that $\widehat{A}$ is locally support-finite. Then $\KG(\widehat{A})\in\{0,2,\infty\}$ and the following assertions hold. 
\begin{enumerate}[\rm(1)]
	\item We have $\KG(\widehat{A})=0$ if and only if $\widehat{A}\cong\widehat{B}$ where $B$ is some tilted algebra of Dynkin type.
	\item We have $\KG(\widehat{A})=2$ if and only if $\widehat{A}\cong\widehat{B}$ where $B$ is some tilted algebra of Euclidean type.
	\item We have $\KG(\widehat{A})=\infty$ if and only if $\widehat{A}$ is wild or $\widehat{A}\cong\widehat{B}$ where $B$ is some tubular algebra.
\end{enumerate} 
\end{thm}

{\bf Proof.} We show all assertions simultaneously. Observe that Theorem 7.1 yields the repetitive category $\widehat{B}$ is locally support-finite in every case considered in assertions $(1)$, $(2)$ and $(3)$ (see also \cite{Sk0}). Moreover, $\widehat{B}_{n}$ is a finite convex subcategory of $\widehat{B}$, $\widehat{B}_{n}\subseteq\widehat{B}_{n+1}$, for any $n\geq 0$, and $\widehat{B}=\bigcup_{n\in\NN}\widehat{B}_{n}$. 

Assume that $B$ is a tilted algebra of Dynkin type. We show that $\KG(\widehat{B})=0$. Indeed, in this case $\widehat{B}$ is locally representation-finite, see \cite{AHR}. Thus $\widehat{B}_{n}$ is representation-finite, so $\KG(\widehat{B}_{n})=0$ by Auslander's result, for any $n\in\NN$. Hence Lemma 4.3 yields $\KG(\widehat{B})=0$.

Assume that $B$ is a tilted algebra of Euclidean type. We show that $\KG(\widehat{B})=2$. It follows from \cite{AsSk1}, \cite{Sk0} and \cite{ANS} that the tilted algebras of Euclidean type and their repetitive categories are cycle-finite of domestic type. This implies that the algebra $\widehat{B}_{n}$ is cycle-finite of domestic type, for any $n\in\NN$. Consequently, Theorem 7.2 yields $\KG(\widehat{B}_{n})=2$, for any $n\in\NN$. Hence Lemma 4.3 yields $\KG(\widehat{B})=2$.

Assume that $B$ is a tubular algebra. It follows from \cite{Ge2} that $\KG(B)=\infty$. Proposition 4.1 implies that $\mod(B)$ is a contravariantly finite subcategory of $\mod(\widehat{B})$. Hence the category $\CF(B)$ is a Serre subcategory of $\CF(\widehat{B})$ and thus $\KG(B)\leq\KG(\widehat{B})$, see Section 3. This yields $\KG(\widehat{B})=\infty$. 

Assume that $\widehat{A}$ is wild. Then there is a finite convex subcategory $B$ of $\widehat{A}$ such that $B$ is wild, see for example \cite{DoSk2}. Consequently, $\KG(B)=\infty$ and so $\KG(\widehat{A})=\infty$ as in the previous case.

To complete the proof, observe that Theorem 7.1 yields $\KG(\widehat{A})\in\{0,2,\infty\}$. It follows easily that all assertions hold. \epv

\section{Krull-Gabriel dimension of selfinjective algebras}

In this section we determine the Krull-Gabriel dimension of standard selfinjective algebras of polynomial growth. We also discuss the existence of super-decomposable pure-injective modules for this class of algebras. The existence of such modules is related to the Krull-Gabriel dimension.

Note that in a recent paper \cite{Bush} M. Bushell gives a description of the \textit{Ziegler spectrum} (the space of isomorphism types of all indecomposable pure-injective modules, see \cite{Zi} for details) of the same class of algebras. Although the concepts of Krull-Gabriel dimension and Ziegler spectrum are related (for example, via the \textit{Cantor-Bendixson rank}, see Sections 7 and 8 in \cite{Zi}), this relation is not discussed in \cite{Bush}. Since the present paper doesn't discuss it as well, we view the results of this section as independent of \cite{Bush}. 

Recall that an algebra $A$ is \textit{standard} if and only if there exists a Galois covering functor $\Gamma\ra A$ such that $\Gamma$ is a \textit{simply connected} locally bounded category, see \cite{AsSk3}, \cite{Sk0} and \cite{Sk4}. 

It follows from \cite{Sk0}, \cite{Sk4} that representation-infinite standard selfinjective algebras of polynomial growth are orbit algebras of the form $\widehat{B}\slash G$ where $B$ is a tilted algebra of Euclidean type or a tubular algebra and $G$ is an infinite cyclic (so torsion-free) admissible group of $K$-linear automorphisms of $\widehat{B}$. Therefore, in order to determine the Krull-Gabriel dimension of these algebras, it suffices to apply Theorem 7.3 and Theorem 6.3. 

The following theorem determines the Krull-Gabriel dimension of standard selfinjective algebras of polynomial growth. This theorem supports Conjecture 1.1 of M. Prest on the finiteness of Krull-Gabriel dimension, see Section 1.

\begin{thm} Assume that $A$ is a standard selfinjective algebra over an algebraically closed field $K$. Then the following assertions hold.
\begin{enumerate}[\rm (1)]
	\item If the algebra $A$ is representation-infinite domestic, then $\KG(A)=2$.
	\item If the algebra $A$ is nondomestic of polynomial growth, then $\KG(A)=\infty$.
\end{enumerate}
\end{thm}

{\bf Proof.} (1) It follows from \cite{Sk0} that $A$ is isomorphic to the orbit algebra $\widehat{B}\slash G$ where $B$ is a tilted algebra of Euclidean type and $G$ is an admissible infinite cyclic group of $K$-linear automorphisms of $\widehat{B}$. The category $\widehat{B}$ is locally support-finite and thus Theorem 6.3 (2) and Theorem 7.3 (2) yield $\KG(A)=\KG(\widehat{B})=2$.

(2) It follows from \cite{Sk0} that $A$ is isomorphic to the orbit algebra $\widehat{B}\slash G$ where $B$ is a tubular algebra and $G$ is an admissible infinite cyclic group of $K$-linear automorphisms of $\widehat{B}$. The category $\widehat{B}$ is locally support-finite and thus Theorem 6.3 (2) and Theorem 7.3 (3) yield $\KG(A)=\KG(\widehat{B})=\infty$. \epv

\begin{cor} Assume that $A$ is a standard representation-infinite selfinjective algebra of polynomial growth over an algebraically closed field $K$. Then $\KG(A)=2$ if and only if the infinite radical $\rad_{A}^{\omega}$ is nilpotent.
\end{cor}

{\bf Proof.} The main theorem of \cite{KerSk} implies that $A$ is domestic if and only if $\rad_{A}^{\omega}$ is nilpotent. Hence the assertion follows directly from Theorem 8.1. \epv

Assume that $\cal{R}$ is a ring with a unit. An $\cal{R}$-module $M$ is \textit{super-decomposable} if and only if $M$ does not have an indecomposable direct summand. For the concept of \textit{pure-injectivity} we refer to \cite{Kiel}, see also \cite{HZ} and \cite[Section 7]{JeLe}.

The problem of the existence of super-decomposable pure-injective $\cal{R}$-modules is studied for the first time in \cite{Zi}. The case when $\cal{R}$ is a finite dimensional algebra over a field is studied, in particular, in \cite{Pr}, \cite{Pu1}, \cite{Pu2}, \cite{Har}, \cite{Pr3}, \cite{KaPa}, \cite{KaPa2}, \cite{Pa}, \cite{GP} and \cite{KaPa3}. It is conjectured (see for example \cite{Pr2}) that if $\cal{R}$ is a finite dimensional algebra over an algebraically closed field, then $\cal{R}$ is of domestic type if and only if there is no super-decomposable pure-injective $\cal{R}$-module. The following theorem supports this conjecture.

\begin{thm} Assume that $A$ is a standard selfinjective algebra over an algebraically closed field. If $A$ is representation-infinite domestic, then there is no super-decomposable pure-injective $A$-module.
\end{thm}

{\bf Proof.} It follows from Theorem 8.1 (1) that $\KG(A)=2$. Since $\KG(A)$ is finite, super-decomposable pure-injective $A$-modules do not exist, see for example \cite{Pr2}. \epv


Assume that $A$ is a standard selfinjective algebra over an algebraically closed field. If $A$ is nondomestic of polynomial growth, then $\KG(A)=\infty$ from Theorem 8.1 (2). It is not known whether this fact implies the existence of a super-decomposable pure-injective $A$-module. 

\section*{Acknowledgements} The research has been supported by the research grant DEC-2011/02/A/ST1/00216 of the National Science Center Poland. 

We are grateful to Stanis{\l}aw Kasjan and Andrzej Skowro\'nski for some helpful comments concerning the paper. We are indebted to the anonymous referee for all remarks. These suggestions allowed us to improve, to a large extent, the initial version of the paper. In particular, we have formulated the results of Section 6 in greater generality and described the Krull-Gabriel dimension of all locally support-finite repetitive $K$-categories in Section 7.


\begin{thebibliography}{}

\normalsize
\baselineskip=17pt

\bibitem{AHR} I. Assem, D. Happel and O. Rold\'an, Representation-finite trivial extension algebras, J. Pure Appl. Algebra 33 (1984), 235--242.

\bibitem{ANS} I. Assem, J. Nehring and A. Skowro\'nski, Domestic trivial extensions of simply connected algebras, Tsukuba J. Math. 13 (1989), 31--72.

\bibitem{AsSiSk} I. Assem, D. Simson and A. Skowro\'nski, {\em Elements of the Representation Theory of Associative Algebras}, Vol. 1, {\em Techniques of representation theory}, London Mathematical Society Student Texts; 65, Cambridge University Press, 2006.

\bibitem{AsSk1} I. Assem and A. Skowro\'nski, Algebras with cycle-finite derived categories, \textit{Math. Ann.} 280, 441--463 (1988).

\bibitem{AsSk2} I. Assem and A. Skowro\'nski, Minimal representation-infinite coil algebras, \textit{Manuscripta Math.} 67, 305--331 (1990). 

\bibitem{AsSk3} I. Assem and A. Skowro\'nski, On some classes of simply connected algebras, \textit{Proc. London Math. Soc} 56, 417--450 (1998).

\bibitem{AsSk4} I. Assem and A. Skowro\'nski, On tame repetitive algebras, \textit{Fund. Math}. 142 (1993), no. 1, 59--84.

\bibitem{Au0} M. Auslander, Functors and morphisms determined by objects. Representation theory of algebras (Proc. Conf., Temple Univ., Philadelphia, Pa., 1976), pp. 1--244. Lecture Notes in Pure Appl. Math., Vol. 37, Dekker, New York, 1978. 

\bibitem{Au} M. Auslander, A functorial approach to representation theory, in: \textit{Representations of Algebras, Lecture Notes in Math.} Vol 944 (1982), 105--179. 

\bibitem{AuRe} M. Auslander and I. Reiten, Applications of contravariantly finite subcategories, \textit{Adv. Math.} 86 (1991), 111--152.



\bibitem{BobKr} G. Bobi\'nski and H. Krause, The Krull-Gabriel dimension of discrete derived categories, \textit{Bull. Sci. Math.} 139 (2015), no. 3, 269--282.

\bibitem{BoGa} K. Bongartz and P. Gabriel, Covering spaces in representation-theory, {\em Invent. Math.} 65 (1981/82) 331--378.



\bibitem{Bush} M. Bushell, Ziegler Spectra of Self Injective Algebras of Polynomial Growth, arXiv:1712.01575.

\bibitem{BuRi} M. C. R. Butler, C. M. Ringel, Auslander-Reiten sequences with few middle terms and applications to string algebras, {\em Comm. Algebra} 15 (1987) 145--179.


\bibitem{DoSk} P. Dowbor and A. Skowro\'nski, Galois coverings of representation-infinite algebras, {\em Comment. Math. Helvetici} 62 (1987) 311--337.

\bibitem{DoSk2} P. Dowbor and A. Skowro\'nski, On the representaion type of locally bounded categories, {\em Tsukuba J. Math.} 10 (1986), no. 1, 63--72. 


\bibitem{Ga0} P. Gabriel, Des categori\'es ab\'eliennes, {\em Bull. Soc. Math. France} 90 (1962), 323--448.

\bibitem{Ga} P. Gabriel, The universal cover of a representation-finite algebra, in: \textit{Representations of Algebras} (Puebla, 1980), \textit{Lecture Notes in Math.}, 903, Springer, 1981, 68--105.

\bibitem{GP} L. Gregory and M. Prest, Representation embeddings, interpretation functors and controlled wild algebras,  \textit{J. Lond. Math. Soc.} (2) 94 (2016), no. 3, 747--766.

\bibitem{Ge1} W. Geigle, The Krull-Gabriel dimension of the representation theory of tame hereditary artin algebras and application to the structure of exact sequences, {\em Manuscr. Math.}, 54 (1985), 83--106.

\bibitem{Ge2} W. Geigle, Krull dimension of Artin algebras, in: {\em Representation Theory I. Finite-dimensional Algebras,} \textit{Lecture Notes in Math.}, No. 1177, Springer-Verlag, Berlin, Heidelberg, New York, 1986, pp. 135--155.



\bibitem{Har} R. Harland, {\em Pure-injective modules over tubular algebras and string algebras}, PhD thesis, University of Manchester, 2011.

\bibitem{HW} D. Hughes and W. Waschb\"usch, Trivial extensions of tilted algebras, \textit{Proc. London Math. Soc.} 46 (1983), 347--364.

\bibitem{HZ} B. Huisgen-Zimmermann, Purity, algebraic compactness, direct sum decompositions, and representation type, in: \textit{Infinite Length Modules} (Bielefeld, 1998), \textit{Trends Math.}, Birkh{\"a}user, Basel, 2000, 331-367.

\bibitem{JeLe} Ch. U. Jensen and H. Lenzing, {\em Model Theoretic Algebra with particular emphasis on Fields, Rings, Modules}, Algebra, Logic and Applications, Vol. 2, Gordon \& Breach Science Publishers, 1989.


\bibitem {KaPa} S. Kasjan and G. Pastuszak, On two tame algebras with super-decomposable pure-injective modules, {\em Colloq. Math.} 123 (2011) 249--276.

\bibitem{KaPa2} S. Kasjan and G. Pastuszak, On the existence of super-decomposable pure-injective modules over strongly simply connected algebras of non-polynomial growth, {\em Colloq. Math.} 136 (2014) 179--220.

\bibitem{KaPa3} S. Kasjan and G. Pastuszak, Super-decomposable pure-injective modules over algebras with strongly simply connected Galois coverings, \textit{J. Pure Appl. Algebra}, Vol. 220 no. 8 (2016) 2985--2999.

\bibitem{KerSk} O. Kerner and A. Skowro\'nski, On module categories with nilpotent infinite radical, \textit{Compositio Math.} 77 (1991), no. 3, 313--333.


\bibitem{Kiel} R. Kie\l pi\'nski, On $\Gamma $-pure-injective modules, \textsl{Bull. Acad. Polon. Sci. S\'er. Sci. Math. Astronom. Phys.} 15  (1967) 127--131.

\bibitem{Kr2} H. Krause, Generic modules over artin algebras, \textit{Proc. London Math. Soc.} 76 (1998), 276--306.

\bibitem{Kr} H. Krause, The spectrum of a module category, \textit{Mem. Amer. Math. Soc.} 149 (2001), no. 707.


\bibitem{LaPrPu} R. Laking, M. Prest and G. Puninski, Krull-Gabriel dimension of domestic string algebras, \textit{Trans. Amer. Math. Soc.}, to appear, arXiv:1506.07703.

\bibitem{McL} S. MacLane,  {\em Categories for the working mathematician}, Graduate Texts in Mathematics, 5. Springer-Verlag, New York, 1998.

\bibitem{Mal} P. Malicki,  Krull dimension of tame generalized multicoil algebras, \textit{Algebr. Represent. Theory} 18 (2015), no. 3, 881--894.

\bibitem{MP} R. Martinez-Villa and J. A. de la Pe\~na, The universal cover of a quiver with relations, {\em J. Pure. Appl. Algebra} 30 (1983). 277--292.





\bibitem{Pa} G. Pastuszak, Strongly simply connected algebras with super-decomposable pure-injective modules, \textit{J. Pure Appl. Algebra} 219 (2015), no. 8, 3314--3321.




\bibitem{Pog} Z. Pogorza\l y, On star-free bound quivers, {\em Bull. Polish Acad. Sci. Math.} 37 (1989), 255--267.

\bibitem{Po} N. Popescu, {\em Abelian Categories with Applications to Rings and Modules}, L.M.S. Monographs, Academic Press, London, 1973.

\bibitem{Pr} M. Prest, {\em Model theory and modules}, London Mathematical Society Lecture Note Series 130, Cambridge University Press, Cambridge, 1988.

\bibitem{Pr2} M. Prest, {\em Purity, Spectra and Localization}, Encyclopedia of Mathematics and Its Applications 121, Cambridge University Press, Cambridge, 2009.

\bibitem{Pr3} M. Prest, Superdecomposable pure-injective modules, in: \textit{Advances in Representation Theory of Algebras}, EMS Series of Congress Reports (2013) 263--296.



\bibitem{PrPu} M. Prest and G. Puninski, Ringel's conjecture for domestic string algebras, \textit{Math. Z.} 282 (2016), no. 1-2, 61--77.

\bibitem{Pu1} G. Puninski, Superdecomposable pure-injective modules exist over some string algebras, {\em Proc. Amer. Math. Soc. } 132 (2004),  1891--1898.

\bibitem{Pu2} G. Puninski, How to construct a 'concrete' superdecomposable pure-injective module over a string algebra, {\em J. Pure Appl. Algebra} 212 (2008), 704--717.

\bibitem{Pu3} G. Puninski, Pure injective indecomposable modules over 1-domestic string algebras, \textit{Algebr. Represent. Theory} 17 (2014), no. 2, 643--673.

\bibitem{PPT} G. Puninski, V. Puninskaya and C. Toffalori,  Krull-Gabriel dimension and the model-theoretic complexity of the category of modules over group rings of finite groups.  {\em J. Lond. Math. Soc.}  78 (2008), 125--142.


\bibitem{Ri} C. M. Ringel, \emph{Tame Algebras and Integral Quadratic Forms}, Lecture Notes in Math., No. 1099, Springer-Verlag, Berlin, Heidelberg, New York, 1984.




\bibitem{Sch1} J. Schr\"oer, On the infinite radical of a module category, \textit{Proc. London Math. Soc.} (3) 81 (2000) 651--674.

\bibitem{Sch2} J. Schr\"oer, On the Krull-Gabriel dimension of an algebra, \textit{Math. Z.} 233, (2000) 287--303.

\bibitem{Sch3} J. Schr\"oer, The Krull-Gabriel dimension of an algebra - open problems and conjectures, in: \textit{Infinite Length Modules}, Trends Math. (2000), 419--424.





\bibitem{SiSk3} D. Simson and A. Skowro\'nski, \emph{Elements of the Representation Theory of Associative Algebras} 3: \emph{Representation-Infinite Tilted Algebras}, London Mathematical Society Student Texts 72, Cambridge University Press, 2007.

\bibitem{SkBC} A. Skowro\'nski, Algebras of polynomial growth, in: \textit{Topics in algebra, Banach Center Publ.} 26, Part 1, PWN, Warsaw, 1990, 535--568.

\bibitem{Sk0} A. Skowro\'nski, Selfinjective algebras of polynomial growth, {\em Math. Ann.} 285 (1989), 177--199.

\bibitem{Sk4}  A. Skowro\'nski,  Selfinjective algebras: finite and tame type, in: {\em Contemp. Math.} 406 (2006) 169--238.


\bibitem{Sk1} A. Skowro\'nski, Simply connected algebras of polynomial growth, {\em Compositio Math.}  109  (1997),  99--133.

\bibitem{Sk3} A. Skowro\'nski, Tame algebras with strongly simply connected Galois coverings, {\em Colloq. Math.} 72 (1997), 335--351.

\bibitem{Sk5} A. Skowro\'nski, The Krull-Gabriel Dimension of Cycle-Finite Artin Algebras, \textit{Algebr. Represent. Theory}, Vol. 19 no. 1, 215--233 (2016). 

\bibitem{SkWa} A. Skowro\'nski and J. Waschb\"usch, Representation-finite biserial algebras, {\em J. Reine Angew. Math.} 345 (1983) 172--181.



\bibitem{We} M. Wenderlich, Krull dimension of strongly simply connected algebras, {\em Bull. Polish Acad. Sci., Ser. Math.}, 44 (1996), 473--480.

\bibitem{Zi} M. Ziegler, Model theory of modules, {\em Annals of Pure and Applied Logic} 26 (1984), 149--213.

\end{thebibliography}
\end{document}